\documentclass[11pt]{article}

\usepackage[margin=1in]{geometry}
\usepackage{amsmath,amssymb,amsthm}
\usepackage{microtype}
\usepackage[shortlabels]{enumitem}
\usepackage[colorlinks=true,linkcolor=blue,citecolor=blue,urlcolor=blue]{hyperref}
\hypersetup{hypertexnames=false,pdftitle={What the occupation law of a finite Markov chain determines},pdfauthor={Arjun Pemmasani}}

\newtheorem{theorem}{Theorem}[section]
\newtheorem{lemma}[theorem]{Lemma}
\newtheorem{proposition}[theorem]{Proposition}
\newtheorem{corollary}[theorem]{Corollary}
\newtheorem{introtheorem}{Theorem}

\theoremstyle{definition}
\newtheorem{definition}[theorem]{Definition}
\newtheorem{example}[theorem]{Example}
\theoremstyle{remark}
\newtheorem{remark}[theorem]{Remark}

\newcommand{\R}{\mathbb{R}}
\newcommand{\E}{\mathbb{E}}

\newcommand{\supp}{\operatorname{supp}}
\newcommand{\ones}{\mathbf{1}}

\title{What the occupation law of a finite Markov chain determines:\\ time reversal and principal minors}
\author{Arjun Pemmasani\\[2pt]
\small Department of Mathematics, Harvey Mudd College\\
\small Claremont, CA 91711, USA\\
\small \texttt{apemmasani@hmc.edu}}
\date{}

\begin{document}
\maketitle

\begin{abstract}
What does the occupation law of a finite Markov chain determine about its generator $Q$? From the
stationary start, the law of the occupation times at one time determines $Q$ up to time
reversal and reversals of memoryless blocks, the sets of states that are all left at one rate with
one landing law and are entered in proportion to the stationary law. These reversals commute and
give exactly $2^k$ generators, with $k$ read off a tree of memoryless blocks. This leaves only time
reversal when $Q$ has no rank-one cut, when there are at most 4 states, and generically. The law
at one time determines the law at every time, and a generic start determines $Q$. A renewal block
is a set of states that is left at the same rates from each of its states, is entered with one
law, and is charged by the start in proportion to that law. From every start the law determines
$Q$ up to lumped reversals, which lump disjoint renewal blocks to single states and reverse
a memoryless block of the lumped chain. We also find every start that hides time reversal. At the
scale of large deviations, the functional of Donsker and Varadhan determines exactly the principal
minors of $-Q$.

\end{abstract}



\section{Introduction}\label{sec:intro}

If we record how long a continuous-time Markov chain on $[d]=\{1,\dots,d\}$ spends in each of its
states, and omit the order of its visits and jumps, what can we determine about the generator $Q$?
Let $O(T)$ be the vector of the times that the chain spends in each state up to
time $T$. We answer the question for the law of $O(T)$, and for its logarithmic scale, the
functional $J_V$ of Donsker and Varadhan~\cite{dv1975}, which is the rate function of $O(T)/T$.

A stationary chain and its time reversal spend the same time in each state. Hence, from the
stationary start, the law of $O(T)$ never separates $Q$ from its time reversal. We show that for
generic $Q$ it separates $Q$ from every other generator. For an elliptic diffusion the analogous
question has a classical answer at the scale of large deviations. Donsker and Varadhan~\cite{dv1976} showed that the principal
eigenvalues of the generator under all potentials, on all subdomains, determine it up to a
ground-state transform and time reversal. In a companion paper on diffusions~\cite{pemmasaniC3}
we show that on a closed manifold of dimension at least 2 a single central limit covariance of the
occupation integrals already leaves only time reversal. For finite chains, Kim and Landim~\cite[Section~5]{kl2024}
recover a reversible $Q$ from $J_V$ (their Corollary~5.3) and show with the uniform 3-cycle that
$J_V$ does not determine a nonreversible one (their Example~5.4). For a finite chain less is
determined than for a diffusion. The functional $J_V$ determines only the principal minors of $-Q$,
and the stationary law of $O(T)$ can fail to separate $Q$ from generators other than its time
reversal.

\begin{table}[ht]
\centering
\small
\begin{tabular}{p{0.36\textwidth}p{0.36\textwidth}l}
\hline
Observation & What it determines & Source\\
\hline
the path & $Q$ & \\
the law of $O(T)$ at one $T$, a start $\mu$ without renewal blocks & $Q$ if $\mu\ne\pi$, and $Q$ up
to time reversal if $\mu=\pi$ & Thm.~\ref{thm:intro-one-time}\\
the law of $O(T)$ at one $T$, stationary start & $Q$ up to $(\mathbb Z/2)^k$, generated by time
reversal and commuting reversals of memoryless blocks & Thm.~\ref{thm:intro-fiber}\\
the law of $O(T)$ at one $T$, any start $\mu$ & $Q$ up to lumped reversals &
Thm.~\ref{thm:intro-one-time}\\
the functional $J_V$ & the principal minors of $-Q$ & Prop.~\ref{prop:intro-dv}\\
\hline
\end{tabular}
\caption{What the observations of a finite irreducible chain determine about its generator.}
\label{tab:hierarchy}
\end{table}

\subsection{The stationary fiber}\label{sec:intro-law}

Write $\pi$ for the stationary law of an irreducible $Q$ and $Q^*=D_\pi^{-1}Q^\top D_\pi$ for its
time reversal, with $D_\pi$ the diagonal matrix of $\pi$. Let $B$ be a set of states with
$2\le|B|\le d-1$. We call $B$ a memoryless block of $Q$ if 2 conditions hold. First, all states of
$B$ jump to each state outside $B$ at the same rate. Second, every jump into $B$ from outside lands
with the law $\pi_B/\pi(B)$. For such a $B$, the reversal of $Q$ inside $B$ replaces the block
$Q_B$ by $D_{\pi_B}^{-1}Q_B^\top D_{\pi_B}$ and keeps the other entries of $Q$. We call a partition
$[d]=X\cup Y$ with $|X|,|Y|\ge2$ a cut of $Q$ if the blocks $Q[X,Y]$ and $Q[Y,X]$ have rank at
most 1.

Two sets overlap if they meet and neither contains the other. The memoryless blocks that overlap
no other memoryless block, together with $[d]$, form a tree under inclusion. We call a node $N$ of
this tree active if $\pi_iq_{ij}\ne\pi_jq_{ji}$ for some states $i$ and $j$ of $N$ that lie together
in no smaller node. Let $k$ be the number of active nodes.

\begin{introtheorem}[the stationary fiber]\label{thm:intro-fiber}
Let $Q$ and $Q'$ be irreducible generators on $[d]$, each started from its stationary law.
\begin{enumerate}[(a)]
\item $O(T)$ has the same law under $Q$ and $Q'$, for one $T>0$ or equivalently for every $T$, if
and only if $Q'$ is reached from $Q$ or from $Q^*$ by reversals of memoryless blocks, each a
memoryless block of the current generator.
\item There are exactly $2^k$ such $Q'$. The reversals commute, and $(\mathbb Z/2)^k$ acts on these
generators simply transitively.
\item They are only $Q$ and $Q^*$ when $d\le4$, when $Q$ has no cut, and for all positive rates
outside a finite union of proper linear subspaces. For $d=5$ there are others.
\end{enumerate}
\end{introtheorem}

This is Theorems~\ref{thm:stationary-fiber} and~\ref{thm:fiber-count}, Corollary~\ref{cor:stationary-cases}
and Example~\ref{ex:five-state}. When $Q$ has no memoryless block, as for generic rates, and when
$d\le3$, the proof needs only the theorems of Loewy and of Hartfiel and Loewy
(Remark~\ref{rem:no-cggr}). Inside a memoryless block each sojourn starts in the stationary law of
the block and leaves it in a way that does not depend on the state, so it can be reversed. For the
converse we border $Q$ by $\pi$. Every cut of the bordered matrix is a memoryless block
(Lemma~\ref{lem:bordered-cuts}), and the theorem of Chatterjee, Ghosh, Gurjar and
Raj~\cite{cggr2024} on matrices with equal principal minors does the rest.

\subsection{One time and every start}\label{sec:intro-onetime}

A renewal block of a start $\mu$ is a set $B$ of at least 2 and at most $d-1$ states that every
state of $B$ leaves at the same rates, that every jump into $B$ enters with one law $\rho$, and that
$\mu$ charges in proportion to $\rho$, possibly not at all. A lumped reversal lumps disjoint renewal
blocks to single states and then reverses a memoryless block of the lumped chain that the lumped
start charges in proportion to its stationary law, or reverses the whole lumped chain when the
lumped start is stationary.

\begin{introtheorem}[one time and every start]\label{thm:intro-one-time}
Let $Q$ be an irreducible generator on $[d]$ and $Q'$ a generator.
\begin{enumerate}[(i)]
\item If $O(T_0)$ has the same law under $Q$ from a start $\mu$ and under $Q'$ from a start
$\mu'$, for one $T_0>0$, then it has the same law for every $T$, $\mu'=\mu$, and $-Q'$ has the
principal minors of $-Q$.
\item If $\mu$ has no renewal block, then the law of $O(T_0)$ at one $T_0>0$ from $\mu$ determines
$Q$, or $Q$ up to time reversal when $\mu=\pi$. This holds for every start outside a finite union of
proper linear subspaces, which depend on $Q$, and for some starts with zeros.
\item For every start $\mu$, the law of $O(T_0)$ at one $T_0>0$ from $\mu$ determines $Q$ exactly
up to finitely many lumped reversals.
\end{enumerate}
\end{introtheorem}

This is Theorems~\ref{thm:one-horizon}, \ref{thm:generic-start} and~\ref{thm:every-start}. The
renewal blocks are the cuts of $Q$ bordered by the column $\ones$ and the row $\mu$. Such a set can hide something, since lumping it can produce a memoryless block
(Example~\ref{ex:lumped}). For the converse in (iii) we show that the minimal renewal blocks, and
those with 2 states, are renewal blocks of every generator with the same law, and we lump one of
them and induct on $d$. With $\mu=\pi$ the renewal blocks are the memoryless blocks, and a lumped
reversal is a composition of reversals of memoryless blocks and of $Q\mapsto Q^*$: it reverses the
set of states that the reversed lumped block stands for, or the whole chain, and then reverses
again each lumped renewal block inside it. So Theorem~\ref{thm:intro-fiber}(a) is the case
$\mu=\pi$ of (iii). We keep its own proof, which is short, rests on the stated theorem of
Chatterjee, Ghosh, Gurjar and Raj~\cite[Theorem~1.1]{cggr2024}, and gives the count $2^k$. The
proof of (iii) rests on more: two lemmas from inside the proof of that theorem~\cite[Lemmas~3.5
and~3.6]{cggr2024}, which are in the arXiv version and in the STOC 2025 version of their paper, and
which we transfer to matrices with zero entries (Lemma~\ref{lem:cggr-transfer}).
Theorems~\ref{thm:intro-fiber}
and~\ref{thm:intro-one-time} hold for the counts of the discrete chain $I+hQ$ for all $N$ at one
fixed $h$.

\subsection{Starts that hide time reversal}\label{sec:intro-hiding}

The stationary start never separates $Q$ from $Q^*$, and a generic start does. Write a start as
$\mu=\pi\circ f$. A circulation component of $Q$ is a connected component of the graph of the edges
$\{i,j\}$ with $\pi_iq_{ij}\ne\pi_jq_{ji}$, and a closed block is a memoryless block with no net
flux across its boundary.

\begin{introtheorem}[which starts hide time reversal]\label{thm:intro-hiding}
Let $Q$ be irreducible. A start $\mu=\pi\circ f$ gives the same law of $O(T)$ under $Q$ and under
$Q^*$ if and only if $f$ is constant on the smallest closed block containing each circulation
component, and equivalently if and only if $\E_\pi[f(\xi_T)-f(\xi_0)\mid O(T)]=0$ for every $T$.
If some circulation component lies in no closed block, only the stationary start does.
\end{introtheorem}

This is Theorem~\ref{thm:hiding}. The difference of the 2 transforms is a bilinear form in the flux
matrix $D_\pi Q-Q^\top D_\pi$ (Lemma~\ref{lem:flux-form}), and an induction along the construction
of the closed block shows that it sees every state the construction adds.

\subsection{The functional of Donsker and Varadhan}\label{sec:intro-dv}

At the logarithmic scale the answer is a statement about principal minors.

\begin{proposition}[what the functional determines]\label{prop:intro-dv}
Let $Q$ and $Q'$ be irreducible generators on $[d]$. Then $J_V=J'_V$ on the simplex if and only if
$-Q$ and $-Q'$ have the same principal minors. In that case the exit rate $\lambda_F(v)$ of every
face $F$ under every field $v$, and the functional $J_F$ of the chain killed on leaving $F$, are
the same for $Q$ and $Q'$.
\end{proposition}

This is Proposition~\ref{prop:dv}, by the variational formula of Donsker and
Varadhan~\cite{dv1975} for the killed chains and an induction over the faces, which is a case of a
theorem of Delmas, Dronnier and Zitt~\cite[Theorem~2.1]{ddz2024}. The generators with the
principal minors of $-Q$ are $Q$ and $Q^*$ when $Q$ has no cut and when $d\le3$
(Corollary~\ref{cor:generator-fiber}), and in general they are described by the theorem of
Chatterjee, Ghosh, Gurjar and Raj~\cite{cggr2024}. By Theorem~\ref{thm:intro-one-time}(i) the law of
$O(T)$ at one time, from any start, determines the principal minors, so it determines $J_V$.

\subsection{Related work}\label{sec:related}

The exact law of the occupation times of a finite chain goes back to Good~\cite{good1961}, and the
joint rate of the empirical measure and the flow is that of Bertini, Faggionato and
Gabrielli~\cite{bfg2014,bfg2015}. The matrix theory behind the fiber is credited in
Section~\ref{sec:fiber}. Cohen~\cite{cohen1981} and Friedland~\cite{friedland1981} study the
convexity of $v\mapsto\varrho(A+\operatorname{diag}v)$. Brydges, van der Hofstad and
K\"onig~\cite[Theorem~2.1]{bhk2007} write the density of the local times at a fixed time through
cofactors of the generator. In~\cite{pemmasaniRS} we study the occupation counts of a chain that
rarely switches states, which is a forward problem. We compute the local laws of the counts, and we
find their support and their boundary layer.

Identifiability from coarse observations has a long history for aggregated Markov chains, where
one observes which class of states the chain occupies at each time, as in models of ion channels.
Fredkin and Rice~\cite{fredkinrice1986}, Ryd\'en~\cite{ryden1996}, Larget~\cite{larget1998} and
Bruno, Yang and Pearson~\cite{brunoyangpearson2005} study when 2 generators give the same law of
the lumped path. Xiang, Zhou, Deng and Yang~\cite{xzdy2024} identify the generator of a stationary
chain from hitting-time distributions observed at some states. An aggregated chain is seen in the
order of the visits that it makes, but with each class's states merged. In our problem, every
state is seen, but only through the total time spent in each state. Hence, neither observation
determines the other.

Section~\ref{sec:minors} proves Proposition~\ref{prop:intro-dv}, Section~\ref{sec:exactlaw} proves
Theorems~\ref{thm:intro-fiber} to~\ref{thm:intro-hiding}, and Section~\ref{sec:verify} lists the
computations and the open questions.

\section{Setting}\label{sec:prelim}

\emph{The chain.} Let $d\ge2$, and let $Q$ be a generator on $[d]=\{1,\dots,d\}$ with rates
$q_{ij}\ge0$ for $i\ne j$ and exit rates $q_i=\sum_{j\ne i}q_{ij}$, so that $Q_{ii}=-q_i$. Let
$q_{\max}=\max_iq_i$. The start $\mu$ is a probability law on $[d]$ and $V=[d]$ is the set of
all the states. If $Q$ is irreducible, $\pi$ is its stationary law. Let $D_\pi$ be the diagonal
matrix of $\pi$, and let $Q^*=D_\pi^{-1}Q^\top D_\pi$ be its time reversal.

\emph{Faces.} A face is a nonempty set $F\subseteq[d]$. Let $\Delta_F$ be the set of probability
vectors on $[d]$ that vanish off $F$, and $\Delta_F^\circ$ the set of those that are positive on
$F$. The killed generator of $F$ is $Q_F=(Q_{ij})_{i,j\in F}$. Since $q_i$ includes
the rates of leaving $F$, $Q_F$ generates the chain killed when it leaves $F$. For $A\subseteq F$,
the principal submatrix of $-Q_F$ on $A$ is $-Q_A$. The graph $G_F$ has an edge $i\to j$ whenever
$i,j\in F$, $i\ne j$ and $q_{ij}>0$. We call $F$ irreducible if $G_F$ is strongly connected, so
every singleton is irreducible.

\emph{Perron roots.} A square matrix $B$ is Metzler if its off-diagonal entries are nonnegative, and its
Perron root $\varrho(B)$ is its largest real eigenvalue. If $B$ is irreducible, $\varrho(B)$ is
simple with positive left and right eigenvectors, and a positive eigenvector of $B$ belongs to
$\varrho(B)$~\cite[Chapter~2]{bermanplemmons1994}. For a face $F$ and a field
$v\in\R^d$ put
\[
 \lambda_F(v)=-\varrho\bigl(Q_F+\operatorname{diag}v_F\bigr),\qquad\lambda_F=\lambda_F(0).
\]
The number $\lambda_F$ is the exit rate of $F$, and $\lambda_{\{i\}}(v)=q_i-v_i$. If $F$ is
reducible, then we order the strongly connected components $C_1,\dots,C_r$ of $G_F$ such that every
edge between components is going forward. Then, $Q_F$ is block triangular and has diagonal blocks
$Q_{C_j}$, so
$\lambda_F(v)=\min_j\lambda_{C_j}(v)$. The shift
rule $\lambda_F(v+c\ones)=\lambda_F(v)-c$ holds for every face. If $F$ is irreducible, let $r$ and
$l$ be positive right and left Perron vectors of $Q_F+\operatorname{diag}v_F$ with $l\cdot r=1$,
and put $\alpha^*_F(v)=l\circ r$, where $x\circ y$ is the entrywise product. Since the Perron root
is simple, $v\mapsto\lambda_F(v)$ is real-analytic on $\R^F$, and first-order perturbation gives
\begin{equation}\label{eq:perron-derivative}
 \frac{\partial\lambda_F}{\partial v_i}(v)=-\alpha^*_F(v)_i<0,\qquad i\in F .
\end{equation}

\begin{lemma}[monotonicity in the face]\label{lem:face-monotone}
If $F\subsetneq G$ are faces and $G$ is irreducible, then $\lambda_G(v)<\lambda_F(v)$ for every
field $v$. In particular every proper face of an irreducible chain has $\lambda_F>\lambda_V=0$.
\end{lemma}

\begin{proof}
Choose $c$ with $B=Q_G+\operatorname{diag}v_G+cI\ge0$. Its principal submatrix on $F$ is
$Q_F+\operatorname{diag}v_F+cI$, since the diagonal of a killed generator keeps the full exit
rates. Let $\widetilde B$ be this submatrix padded with zeros to $G\times G$. Then
$0\le\widetilde B\le B$, and $\widetilde B\ne B$ because the irreducible matrix $B$ has a positive
entry in a row of $F$ and a column off $F$. So
$\varrho(\widetilde B)<\varrho(B)$~\cite[Chapter~2]{bermanplemmons1994}, and $\varrho(\widetilde B)$
is the Perron root of the submatrix on $F$, since that root is nonnegative and the padding adds
only the eigenvalue $0$. This holds whether or not $F$ is irreducible, and subtracting $c$ gives
the first claim. For $G=V$ and $v=0$,
the vector $\ones$ is a positive null vector of $Q$, so $\lambda_V=0$.
\end{proof}

\emph{The functional.} For $\alpha\in\Delta_F$ let
\[
 J_F(\alpha)=\sup_{u\in(0,\infty)^F}\sum_{i\in F}\alpha_i\Bigl(-\frac{(Q_Fu)_i}{u_i}\Bigr),
\]
the Donsker and Varadhan functional of the chain killed on leaving $F$~\cite{dv1975}, and write
$J_V$ for that of the whole chain. It is the rate function of the occupation measure
$T^{-1}O(T)$ of Section~\ref{sec:exactlaw} as $T\to\infty$~\cite{dv1975}. We use 2 facts.

\begin{proposition}[the variational formula]\label{prop:variational}
Let $F$ be a face.
\begin{enumerate}[(a)]
\item $J_G=J_F$ on $\Delta_F$ for every face $G\supseteq F$. In particular $J_F=J_V$ on $\Delta_F$.
\item For every $v\in\R^F$
\[
 \lambda_F(v)=\min_{\alpha\in\Delta_F}\bigl(J_F(\alpha)-v\cdot\alpha\bigr).
\]
So $-\lambda_F$ is the Legendre transform of $J_F$, and
$J_F(\alpha)=\sup_{v\in\R^F}\bigl(v\cdot\alpha+\lambda_F(v)\bigr)$ for every $\alpha\in\Delta_F$.
\end{enumerate}
\end{proposition}

\begin{proof}
(a) Let $\alpha\in\Delta_F$ and $u\in(0,\infty)^G$. For $i\in F$,
$-(Q_Gu)_i=q_iu_i-\sum_{j\in F\setminus\{i\}}q_{ij}u_j-\sum_{j\in G\setminus F}q_{ij}u_j\le-(Q_Fu_F)_i$,
so $J_G(\alpha)\le J_F(\alpha)$. Letting $u_j\downarrow0$ for $j\in G\setminus F$ gives the reverse
inequality.

(b) Write $Q_F=L_F-\operatorname{diag}k_F$, where $L_F$ is the generator on $F$ with the rates of $Q$
inside $F$ and $k_i=\sum_{j\notin F}q_{ij}$ is the rate of killing at $i$. Then
$Q_F+\operatorname{diag}v_F=L_F+\operatorname{diag}(v_F-k_F)$, and $J_F(\alpha)$ is the functional of
$L_F$ plus $k\cdot\alpha$. If $F$ is irreducible, the formula is that of Donsker and
Varadhan~\cite{dv1975} for the generator $L_F$ and the potential $v_F-k_F$. In general put
$B=Q_F+\operatorname{diag}v_F$, so that $\varrho(B)=-\lambda_F(v)$. Let $C$ be a strongly connected
component of $G_F$ with $\varrho(B_C)=\varrho(B)$. By the irreducible case and (a),
$\varrho(B)=\max_{\alpha\in\Delta_C}(v\cdot\alpha-J_F(\alpha))$, so $\varrho(B)$ is at most the
maximum over $\Delta_F$. Conversely, for every $\varepsilon>0$ there is $u\in(0,\infty)^F$ with
$Bu\le(\varrho(B)+\varepsilon)u$ entrywise. This is the formula of Collatz and Wielandt, which
holds for every Metzler matrix, since a positive matrix $B+sI+\delta\ones\ones^\top$ has a positive
Perron vector and its Perron root tends to $\varrho(B)+s$ as $\delta\downarrow0$. For
$\alpha\in\Delta_F$ the vector $u$ gives
$J_F(\alpha)\ge\sum_i\alpha_i\bigl(v_i-(Bu)_i/u_i\bigr)\ge v\cdot\alpha-\varrho(B)-\varepsilon$. So
the maximum over $\Delta_F$ is at most $\varrho(B)$, and it is attained since $J_F$ is lower
semicontinuous and $\Delta_F$ is compact.

The function $J_F$ is a supremum of affine
functions of $\alpha$, so it is convex and lower semicontinuous, and it is finite on $\Delta_F$,
since $u=\ones$ shows $J_F\ge0$ and every term is at most $q_{\max}$. Extend it by $+\infty$ off
$\Delta_F$. Then the first sentence says that the convex conjugate of $J_F$ is $-\lambda_F$, and
the theorem of Fenchel and Moreau gives $J_F=J_F^{**}$, which is the second.
\end{proof}

\section{The Donsker and Varadhan functional and the principal minors}\label{sec:minors}

This section proves Proposition~\ref{prop:intro-dv} and describes the generators with given
principal minors. Section~\ref{sec:perron-branch} shows that the Perron root of a face on an open
set of fields determines the principal minors of the face. Section~\ref{sec:dv-minors} deduces
that $J_V$ determines exactly the principal minors of $-Q$, and Section~\ref{sec:fiber} describes
the generators with given minors.

\subsection{The polynomial of a face and the Perron branch}\label{sec:perron-branch}

The principal minors of $-Q_F$ are the numbers $\det(-Q_A)$ for $A\subseteq F$, with
$\det(-Q_\emptyset)=1$. They are the coefficients of the multiaffine polynomial
\[
 R_F(w)=\det\bigl(-Q_F-\operatorname{diag}w\bigr)
 =\sum_{A\subseteq F}\det(-Q_A)\prod_{i\in F\setminus A}(-w_i),\qquad w\in\R^F .
\]
For every field $v$, $\lambda_F(v)$ is a real root of $t\mapsto R_F(v+t\ones)$, and if
$F$ is irreducible it is the smallest real root, since the other real eigenvalues of
$-Q_F-\operatorname{diag}v_F$ exceed the Perron eigenvalue.

\begin{lemma}[irreducibility]\label{lem:RF-irreducible}
If $F$ is irreducible, then $R_F$ is irreducible over $\mathbb C$. If $F$ is reducible,
$R_F$ is the product of the polynomials of the strongly connected components of $G_F$.
\end{lemma}

\begin{proof}
Al Ahmadieh~\cite[Lemma~5.1]{alahmadieh2023} proves, over any field and for every
order, that a square matrix $A$ is irreducible exactly when
$\det(\operatorname{diag}x+A)$ is an irreducible polynomial in $x$. Take the field
$\mathbb C$ and $A=-Q_F$. Then $R_F(w)=\det(\operatorname{diag}(-w)+A)$, and the linear
change of variables $x=-w$ preserves irreducibility. For a reducible face, $-Q_F$ is
block triangular with the blocks of the components on the diagonal, and the determinant
factors accordingly.
\end{proof}

The proof of Al Ahmadieh rests on a theorem of Hartfiel and
Loewy~\cite[Theorem~1]{hartfielloewy1984}, that every entry of the adjugate of
$\operatorname{diag}x+A$ is a nonzero polynomial when $A$ is irreducible. The next lemma turns
the vanishing of a polynomial along one root of an irreducible polynomial into a divisibility.

\begin{lemma}[vanishing on a graph]\label{lem:graph}
Let $\mathcal P\in\mathbb C[s,w]$, with $w=(w_1,\dots,w_n)$, be irreducible and monic in $s$ of
degree $m\ge1$. Let $U\subseteq\R^n$ be nonempty and open, and let $\theta$ be a function on $U$
with $\mathcal P(\theta(w),w)=0$. If $\mathcal R\in\mathbb C[s,w]$ and $\mathcal R(\theta(w),w)=0$
for all $w\in U$, then $\mathcal P$ divides $\mathcal R$.
\end{lemma}

\begin{proof}
Divide by the monic $\mathcal P$ in $s$, so $\mathcal R=\mathcal A\mathcal P+\mathcal B$ with
$\deg_s\mathcal B<m$ and $\mathcal B(\theta(w),w)=0$ on $U$. Suppose $\mathcal B\ne0$. Then
$\mathcal P$ does not divide $\mathcal B$, so by Gauss's lemma they have no common factor over
$\mathbb C(w)$, and the resultant $\rho=\operatorname{Res}_s(\mathcal P,\mathcal B)\in\mathbb C[w]$
is not zero. Since $\mathcal P$ is monic, $\rho(w)$ is the product of $\mathcal B(\sigma,w)$ over
the roots $\sigma$ of $\mathcal P(\cdot,w)$, and one of them is $\theta(w)$. So $\rho$ vanishes on
$U$, which is impossible for a nonzero polynomial.
\end{proof}

\begin{lemma}[the Perron branch]\label{lem:perron-branch}
Let $F$ be irreducible, let $Q'$ be another generator whose face $F$ is irreducible,
and let $c\in\R$. Suppose that $\lambda'_F(v)=\lambda_F(v)+c$ for all $v$ in a nonempty
open set $\mathcal U\subseteq\R^F$. Then $-Q'_F$ and $-Q_F+cI$ have the same principal
minors, and $\lambda'_F=\lambda_F+c$ on all of $\R^F$. In particular, for $c=0$,
$-Q'_F$ and $-Q_F$ have the same principal minors.
\end{lemma}

\begin{proof}
Let $k=|F|$ and put $\widetilde R(w)=\det(-Q_F+cI-\operatorname{diag}w)=R_F(w-c\ones)$, which is
irreducible over $\mathbb C$ by Lemma~\ref{lem:RF-irreducible}. In the variables $(t,v)$ put
$\mathcal P(t,v)=(-1)^k\widetilde R(v+t\ones)$ and $\mathcal P'(t,v)=(-1)^kR'_F(v+t\ones)$. Both are
monic of degree $k$ in $t$, and $\mathcal P$ is irreducible, since it is $\pm\widetilde R$ composed
with the linear surjection $(t,v)\mapsto v+t\ones$. For $v\in\mathcal U$ the number
$\theta(v)=\lambda'_F(v)=\lambda_F(v)+c$ is an eigenvalue of $-Q'_F-\operatorname{diag}v$ and of
$-Q_F+cI-\operatorname{diag}v$, so $\mathcal P(\theta(v),v)=0=\mathcal P'(\theta(v),v)$.
Lemma~\ref{lem:graph} gives that $\mathcal P$ divides $\mathcal P'$, and both are monic of degree
$k$ in $t$, so $\mathcal P'=\mathcal P$. At $t=0$ this says $R'_F=\widetilde R$. Finally
$\lambda'_F(v)$ is the smallest real root of $t\mapsto R'_F(v+t\ones)=\widetilde R(v+t\ones)$, which
is $\lambda_F(v)+c$.
\end{proof}

Cohen~\cite{cohen1981} proved that $v\mapsto\varrho(A+\operatorname{diag}v)$ is convex for a
Metzler matrix $A$, and Friedland~\cite[Theorem~4.1]{friedland1981} showed that for irreducible $A$
it is affine only along $\ones$. Lemma~\ref{lem:perron-branch} goes the other way. It recovers
the principal minors of $A$ from this map on an open set.

\subsection{The functional determines exactly the principal minors}\label{sec:dv-minors}

\begin{proposition}[the functional and the minors]\label{prop:dv}
For two irreducible generators $Q$ and $Q'$ on $[d]$ the following are equivalent.
\begin{enumerate}[(i)]
\item $J_V=J'_V$ on $\Delta_V$.
\item $\lambda_V=\lambda'_V$ on a nonempty open set of fields.
\item $-Q$ and $-Q'$ have the same principal minors.
\end{enumerate}
If they hold, then $\lambda_F=\lambda'_F$ on $\R^F$ and $J_F=J'_F$ on $\Delta_F$ for every
face $F$, and $F$ is irreducible for $Q$ exactly when it is irreducible for $Q'$.
\end{proposition}

\begin{proof}
(i) implies (iii). By Proposition~\ref{prop:variational}(a), $J_F=J_V=J'_V=J'_F$ on $\Delta_F$ for
every face $F$, so Proposition~\ref{prop:variational}(b) gives $\lambda_F=\lambda'_F$, and in
particular $\lambda_F(0)=\lambda'_F(0)$. We show $\det(-Q_F)=\det(-Q'_F)$ by induction on $|F|$,
starting from $\det(-Q_\emptyset)=1$. If $-Q_F$ and $-Q'_F$ have the same proper principal
minors, then $R'_F-R_F$ is the constant $\det(-Q'_F)-\det(-Q_F)$. The number $\lambda_F(0)$ is an
eigenvalue of $-Q_F$, and $\lambda'_F(0)=\lambda_F(0)$ is one of $-Q'_F$, so both polynomials vanish
at $\lambda_F(0)\ones$, and the constant is $0$. With $M=Q+cI$ for $c\ge\max(q_{\max},q'_{\max})$,
this is the case of the theorem of Delmas, Dronnier and Zitt~\cite[Theorem~2.1]{ddz2024} that
nonnegative matrices $M$ and $M'$ with $\varrho(M_F)=\varrho(M'_F)$ for every $F$ have the same
principal minors.

(ii) implies (iii). This is Lemma~\ref{lem:perron-branch} with $F=V$ and $c=0$.

(iii) implies (i) and (ii). Equal principal minors give $R'_F=R_F$ for every face $F$, since the
minors of $-Q_F$ are among those of $-Q$. By Lemma~\ref{lem:RF-irreducible}, $F$ is
irreducible for $Q$ exactly when $R_F$ is irreducible, which is the last claim. Every real
eigenvalue of the Metzler matrix $Q_F+\operatorname{diag}v_F$ is at most its Perron root, so for
every face $\lambda_F(v)$ is the smallest real root of $t\mapsto R_F(v+t\ones)$, and likewise for
$Q'$. So $\lambda'_F=\lambda_F$, which gives (ii). By Proposition~\ref{prop:variational}(b),
$J_F(\alpha)=\sup_v(v\cdot\alpha+\lambda_F(v))$ on $\Delta_F$, and likewise for $Q'$, so
$J'_F=J_F$ on $\Delta_F$. With $F=V$ this is (i).
\end{proof}

Lemma~\ref{lem:perron-branch} is needed for (ii), when only the Perron root of the whole space is
known, and only on an open set of fields. That is also the situation of the law of $O(T)$ at one
time, and Step~5 of the proof of Theorem~\ref{thm:one-horizon} runs the same argument.

\subsection{The generators with given minors}\label{sec:fiber}

Without reversibility the exit rates and the edge products $q_{ij}q_{ji}$ do not determine
$Q$, since $Q^*$ has the same exit rates, the same products and the same stationary law.
The next corollary describes the generators with the principal minors of $-Q$ when $Q$ has no cut.
A matrix $C$ is a Z-matrix if $-C$ is Metzler.

\begin{lemma}[the generator point]\label{lem:generator-point}
Let $Q$ be an irreducible generator and $C$ an irreducible Z-matrix with the principal minors of
$-Q$. Then $-C$ has the spectrum of $Q$, and there is a positive diagonal matrix $D$, unique up to a
factor, such that $-D^{-1}CD$ is a generator. It is irreducible, and every real invertible diagonal
$D'$ for which $-D'^{-1}CD'$ has zero row sums is a multiple of $D$. For $C=-Q$ this generator is
$Q$, and for $C=-Q^\top$ it is $Q^*$.
\end{lemma}

\begin{proof}
The characteristic polynomial
$\det(tI-C)=\sum_{s=0}^d(-1)^st^{d-s}E_s(C)$, with $E_s$ the sum of the principal minors of
order $s$, is that of $-Q$. So $-C$ has the spectrum of $Q$, whose eigenvalue of largest real
part is the simple eigenvalue $0$. The theorem of Perron and Frobenius for the irreducible
Metzler matrix $-C$ gives a vector $x>0$, unique up to a factor, with $Cx=0$. With
$D=\operatorname{diag}x$ the matrix $-D^{-1}CD$ has off-diagonal entries $-c_{ij}x_j/x_i\ge0$
and zero row sums, so it is an irreducible generator. If $D'$ is real, diagonal and
invertible and $-D'^{-1}CD'$ has zero row sums, then $C(D'\ones)=0$, so $D'$ is a multiple of
$D$ and gives the same matrix. For $C=-Q$ the generator is $Q$, and for $C=-Q^\top$, whose null
vector is $\pi$, it is $Q^*$.
\end{proof}

\begin{definition}[generator point]\label{def:generator-point}
The generator $-D^{-1}CD$ of Lemma~\ref{lem:generator-point} is the \emph{generator point} of $C$.
\end{definition}

\begin{corollary}[the generator fiber]\label{cor:generator-fiber}
Let $Q$ be an irreducible generator.
\begin{enumerate}[(a)]
\item Every generator $Q'$ with the principal minors of $-Q$ is irreducible and has the
spectrum of $Q$. The set of such generators contains $Q$ and its time reversal $Q^*$.
\item If $Q$ has no cut, that is, for every partition $V=X\cup Y$ with $|X|,|Y|\ge2$ one of
the blocks $Q[X,Y]$, $Q[Y,X]$ has rank at least $2$, then the set in \textup{(a)} is
$\{Q,Q^*\}$. It is a single point exactly when $Q$ is reversible.
\item If $d\le3$, the set in \textup{(a)} is $\{Q,Q^*\}$ with no condition.
\end{enumerate}
\end{corollary}

\emph{What is borrowed.} The matrix theorems used are Loewy~\cite[Theorem~1]{loewy1986} for
order $n\ge4$ without a cut and Hartfiel and Loewy~\cite[Theorem~3]{hartfielloewy1984} for $n\le3$,
combined in~\cite[Lemma~2.8]{cggr2024}. Chatterjee, Ghosh, Gurjar and
Raj~\cite[Theorem~1.1]{cggr2024} describe the whole fiber by cut-transposes, and
Section~\ref{sec:law-cuts} uses their theorem for the stationary law. What we add here is
Lemma~\ref{lem:generator-point} and the irreducibility of every member.

\begin{proof}
We first record a fact used in (b) and (c).

\emph{(R) Real similarities.} If $-Q'=E(-Q)E^{-1}$ or $-Q'=E(-Q)^\top E^{-1}$ for a complex
invertible diagonal $E$, then $E$ may be taken real. In the first case
$q'_{ij}=e_iq_{ij}/e_j$, so $e_i/e_j=q'_{ij}/q_{ij}$ is real on every edge of the graph of
$Q$. That graph is strongly connected, so every ratio $e_i/e_1$ is a product of such ratios
and is real, and $E/e_1$ is a real diagonal similarity. The second case is the same with
$q_{ji}$ in place of $q_{ij}$.

(a) Let $Q'$ be a generator with the principal minors of $-Q$, and suppose that $-Q'$ is
reducible. Then there is a proper nonempty $A\subset V$ with $Q'[V\setminus A,A]=0$, and
$\det(-Q')=\det(-Q'_A)\det(-Q'_{V\setminus A})$. Both factors are proper principal minors,
so they equal $\det(-Q_A)$ and $\det(-Q_{V\setminus A})$. By Lemma~\ref{lem:face-monotone}
every irreducible component of a proper face has a positive exit rate, so every eigenvalue
of $-Q_A$ has positive real part and $\det(-Q_A)>0$, and likewise
$\det(-Q_{V\setminus A})>0$. This contradicts $\det(-Q')=\det(-Q)=0$. So $Q'$ is irreducible,
and the same argument applies to every real matrix with the principal minors of $-Q$. This also
follows from Lemma~\ref{lem:RF-irreducible}, since such a matrix has the polynomial $R_V$ of
$-Q$. The spectra agree by Lemma~\ref{lem:generator-point}. The time reversal $Q^*$ is a
generator, since $Q^*\ones=D_\pi^{-1}Q^\top\pi=0$ and its off-diagonal entries are
$\pi_jq_{ji}/\pi_i\ge0$, and $-Q^*=D_\pi^{-1}(-Q)^\top D_\pi$ has the principal minors of $-Q$,
since transposition and diagonal similarity preserve them.

(b) For $d\le3$ the claim is (c). Let $d\ge4$. Loewy's theorem~\cite[Theorem~1]{loewy1986}
applies to the irreducible matrix $-Q$ without a cut, with no condition on $-Q'$. It gives
an invertible diagonal $E$ with $-Q'=E(-Q)E^{-1}$ or $-Q'^\top=E(-Q)E^{-1}$, and by (R) we may
take $E$ real. In the first case Lemma~\ref{lem:generator-point} with $C=-Q$ gives $Q'=Q$, and in
the second case $Q'=E^{-1}Q^\top E$ and Lemma~\ref{lem:generator-point} with $C=-Q^\top$ gives
$Q'=Q^*$. Finally $Q^*=Q$ exactly when $\pi_iq_{ij}=\pi_jq_{ji}$ for all $i,j$, that is when $Q$
is reversible.

(c) For $d=2$ a generator is fixed by $q_1$ and $q_2$, which are principal minors. Let
$d=3$ and let $Q'$ be in the fiber, so $Q'$ is irreducible by (a). The theorem of Hartfiel
and Loewy~\cite[Theorem~3]{hartfielloewy1984}, as restated in~\cite[Lemma~2.8]{cggr2024},
says that for an irreducible matrix $A$ of order $2$ or $3$ over a field, every matrix with
the principal minors of $A$ is diagonally similar to $A$ or to $A^\top$ over that field. So
$-Q'=E(-Q)E^{-1}$ or $-Q'=E(-Q)^\top E^{-1}$ with $E$ diagonal and invertible, real by (R),
and Lemma~\ref{lem:generator-point} gives $Q'=Q$ or $Q'=Q^*$, as in (b).
\end{proof}

\section{The exact occupation law}\label{sec:exactlaw}

\subsection{Reductions}\label{sec:law-reductions}

Let $(\xi_t)_{t\ge0}$ be the chain with generator $Q$ and start $\mu$, and let
$O(T)\in[0,T]^d$ be its vector of occupation times,
$O_i(T)=\int_0^T\ones\{\xi_t=i\}\,dt$. For $x\in\mathbb C^d$ put
\[
 F^\mu_Q(x)=\mu\bigl(\operatorname{diag}x-Q\bigr)^{-1}\ones ,\qquad F_Q=F^\pi_Q
 \quad\text{when $Q$ is irreducible,}
\]
a rational function of $x$. We write $W=\operatorname{diag}w$ for $w\in\R^d$ and $D=D_\pi$ for the
diagonal matrix of $\pi$.

\begin{lemma}[reductions]\label{lem:reductions}
Let $Q$ be a generator and $\mu$ a start.
\begin{enumerate}[(a)]
\item $\E_\mu e^{-w\cdot O(T)}=\mu e^{T(Q-W)}\ones$ for $w\in\R^d$. The laws of $O(T)$ for all
$T\ge0$ determine $F^\mu_Q$, and $F^\mu_Q$ determines them.
\item Let $h>0$ with $hq_{\max}\le1$, and let $n(N)$ be the composition of $X_1,\dots,X_N$, where
$X_1\sim\mu$ and the transition matrix is $I+hQ$. For small $|y|$ and nonzero $z\in\mathbb C^d$,
\[
 \sum_{N\ge1}y^N\,\E_\mu\prod_iz_i^{n_i(N)}=h^{-1}F^\mu_Q(x),\qquad x_i=\frac{1/(yz_i)-1}h .
\]
So at a fixed $h$ the laws of $n(N)$ for all $N\ge1$ determine $F^\mu_Q$, and $F^\mu_Q$
determines them.
\end{enumerate}
\end{lemma}

\begin{proof}
(a) The function $u(T)=\E_\cdot\,e^{-w\cdot O(T)}$ solves $u'=(Q-W)u$ with $u(0)=\ones$, by the
Feynman and Kac formula, which for a finite chain is the backward equation. For real $s$ larger
than the real parts of the spectrum of $Q-W$,
$\int_0^\infty e^{-sT}\mu e^{T(Q-W)}\ones\,dT=\mu(sI-Q+W)^{-1}\ones=F^\mu_Q(s\ones+w)$. The law of
$O(T)$ lives on $[0,T]^d$, so it is determined by its Laplace transform on $\R^d$. Hence the laws
for all $T$ give the continuous functions $T\mapsto\mu e^{T(Q-W)}\ones$, these give $F^\mu_Q$ on
a nonempty open subset of $\R^d$, and a rational function is determined there. Conversely, for each
$w$ the function $s\mapsto F^\mu_Q(s\ones+w)$ for large $s$ is the Laplace transform of
$T\mapsto\mu e^{T(Q-W)}\ones$, which it determines.

(b) With $Z=\operatorname{diag}z$ and $P=I+hQ$ we have $\E_\mu\prod_iz_i^{n_i(N)}=\mu Z(PZ)^{N-1}\ones$.
Summing the geometric series gives $\mu\,yZ(I-yPZ)^{-1}\ones=\mu\bigl((yZ)^{-1}-P\bigr)^{-1}\ones$,
and $(yZ)^{-1}-I-hQ=h(\operatorname{diag}x-Q)$. The laws for all $N$ give the power series, hence
$F^\mu_Q$ on an open set of $x$, since $z\mapsto x$ is open for fixed $y$. Conversely $F^\mu_Q$ gives
the power series.
\end{proof}

By (a) and (b), each result of this section that assumes equal laws of $O(T)$ for all $T$, or at
one $T$ (which is the same by Theorem~\ref{thm:one-horizon}), holds verbatim for equal laws of
the counts of $I+hQ$ for all $N$ at one fixed $h$.

\subsection{One time determines every time}\label{sec:law-one-time}

Put $\mathcal P(s,w)=\det(sI-Q+W)$. In the variables $(s,s\ones+w)$ it is the polynomial
$R_V(-(s\ones+w))$ of Section~\ref{sec:perron-branch}, which does not involve $s$. So if $Q$ is
irreducible, $\mathcal P$ is irreducible over $\mathbb C$ by Lemma~\ref{lem:RF-irreducible}, and it
is monic of degree $d$ in $s$, so Lemma~\ref{lem:graph} applies to it.

\begin{theorem}[one time]\label{thm:one-horizon}
Let $Q$ be an irreducible generator and $Q'$ any generator on $[d]$, with starts $\mu$ and $\mu'$.
If $O(T_0)$ has the same law under $(Q,\mu)$ and under $(Q',\mu')$ for one $T_0>0$, then $-Q'$
has the principal minors of $-Q$, so $Q'$ is irreducible, and $F^{\mu'}_{Q'}=F^\mu_Q$. So the laws
agree at every $T$.
\end{theorem}

\begin{proof}
\emph{Step 1, the data.} By Lemma~\ref{lem:reductions}(a) the hypothesis says
$\mu e^{T_0(Q-W)}\ones=\mu'e^{T_0(Q'-W)}\ones$ for all $w\in\R^d$.

\emph{Step 2, one exponential.} Let $q=\max(q_{\max},q'_{\max})>0$, and for $v\in\R^d$ with
distinct entries let $\delta(v)$ be the least distance between 2 of them. By the theorem of
Gershgorin the eigenvalues of $V-\epsilon Q$, $V=\operatorname{diag}v$, lie within $2|\epsilon|q$
of the entries of $v$. So for $|\epsilon|<\epsilon_0(v)=\delta(v)/(4q)$ they are simple, given by
analytic functions $\kappa_i(\epsilon)$ with $\kappa_i(0)=v_i$, and so are their spectral
projectors $E_i(\epsilon)$. The same holds for $Q'$. Fix $v^0$ with distinct entries and $v^0_1$
the smallest, and a ball $\mathcal B$ around $v^0$ whose closure consists of such vectors. Then
$\epsilon_*=\min\epsilon_0>0$ on $\mathcal B$, because $\delta$ is continuous and positive on the
closure. This is what makes the set $U$ below open.

For $v\in\mathcal B$ and $t>1/\epsilon_*$ the eigenvalues of $Q-tV$ are $-t\kappa_i(1/t)$ and
\[
 e^{T_0tv_1}\mu e^{T_0(Q-tV)}\ones=\phi_v(1/t)+O(e^{-T_0t\delta(v)/2}),\qquad
 \phi_v(\epsilon)=\mu E_1(\epsilon)\ones\,\exp\Bigl(-T_0\,\frac{\kappa_1(\epsilon)-v_1}{\epsilon}\Bigr),
\]
where $\phi_v$ is analytic for $|\epsilon|<\epsilon_*$. With $\phi'_v$ formed from $Q'$ and
$\mu'$, Step 1 gives $\phi_v(\epsilon)-\phi'_v(\epsilon)=O(e^{-c/\epsilon})$ as $\epsilon\downarrow0$.
Every Taylor coefficient of $\phi_v-\phi'_v$ at $0$ is therefore $0$, and $\phi_v=\phi'_v$. On the
open cone $U=\{tv\mid v\in\mathcal B,\ t>1/\epsilon_*\}$ let $\theta(w)$ be the eigenvalue of $Q-W$
within $2q$ of $-w_1$, which is simple and real since its conjugate is also an eigenvalue there,
and let $\chi(w)=\mu E(w)\ones$ with $E(w)$ its spectral projector. Define $\theta'$ and $\chi'$
from $Q'$ and $\mu'$. We have shown that
\begin{equation}\label{eq:one-exponential}
 \chi(w)\,e^{T_0\theta(w)}=\chi'(w)\,e^{T_0\theta'(w)},\qquad w\in U .
\end{equation}

\emph{Step 3, the residue.} Write $\mu(sI-Q+W)^{-1}\ones=\mathcal N_\mu/\mathcal P$ with
$\mathcal N_\mu=\mu\operatorname{adj}(sI-Q+W)\ones$. Then
$\chi=\mathcal N_\mu(\theta,w)/\partial_s\mathcal P(\theta,w)$, the residue at the simple pole
$s=\theta(w)$. The polynomial $\mathcal N_\mu$ has degree $d-1$ in $s$ and leading coefficient
$\mu\cdot\ones=1$, so $\mathcal P$ does not divide it, and by Lemma~\ref{lem:graph}
$\mathcal N_\mu(\theta(w),w)$ does not vanish identically on any open subset of $U$. So there is a
ball $U_1\subseteq U$ on which $\chi\ne0$, and then $\chi'\ne0$ by~\eqref{eq:one-exponential}.

\emph{Step 4, transcendence.} Fix $w^0\in U_1$ and $\zeta\in\R^d$, and put
$f(r)=T_0(\theta'-\theta)(w^0+r\zeta)$ and $g(r)=(\chi/\chi')(w^0+r\zeta)$ for real $r$ near $0$.
Both are real-analytic, and both are algebraic over $\mathbb C(r)$, since $\theta$ and $\theta'$
are roots of polynomials in $s$ with coefficients in $\mathbb C[r]$, and $\chi$, $\chi'$ are
rational in them and $r$. By~\eqref{eq:one-exponential}, $e^f=g$. We claim that $f$ is constant.
The field $\mathbb C(r,f,g)$ is the field of meromorphic functions of a compact Riemann surface,
on which $f$ and $g$ are meromorphic, and $e^f=g$ continues to the surface minus the poles of $f$
and $g$. If $f$ were not constant it would have a pole $p$, with $f=a\tau^{-k}(1+O(\tau))$, $a\ne0$,
$k\ge1$, in a local coordinate $\tau$ at $p$. Along the directions where $a\tau^{-k}>0$ the
function $|e^f|$ grows faster than every power of $1/|\tau|$, while $g$ is meromorphic at $p$.
So $f$ is constant. As $\zeta$ is arbitrary, the gradient of $\theta'-\theta$ vanishes on $U_1$,
and $\theta'=\theta+c$ there for a real constant $c$.

\emph{Step 5, the minors.} Put $\mathcal P'(s,w)=\det(sI-Q'+W)$. The polynomial
$\mathcal P'(s+c,w)$ vanishes at $s=\theta(w)$ for every $w\in U_1$. This holds only for real $w$,
and Lemma~\ref{lem:graph} turns it into an identity of polynomials. It gives that $\mathcal P$
divides $\mathcal P'(s+c,w)$, and both are monic of degree $d$ in $s$, so they are equal. At $s=0$
this says that $-Q'+cI$ and $-Q$ have the same principal minors, and at $w=0$ that the spectrum
of $Q'$ is that of $Q$ shifted by $c$. Every eigenvalue of a generator has real part at most $0$.
Since $0$ is an eigenvalue of $Q'$, $-c$ is one of $Q$, so $c\ge0$, and since $0$ is an eigenvalue
of $Q$, $c$ is one of $Q'$, so $c\le0$. Hence $c=0$, $-Q'$ has the principal minors of $-Q$, and
$Q'$ is irreducible by Corollary~\ref{cor:generator-fiber}(a).

\emph{Step 6, the transforms.} Now $\mathcal P'=\mathcal P$, $\theta'=\theta$, and $\chi=\chi'$ on
$U_1$ by~\eqref{eq:one-exponential}. So $\mathcal N_\mu-\mathcal N'_{\mu'}$ vanishes at
$s=\theta(w)$ on $U_1$, and by Lemma~\ref{lem:graph} it is divisible by $\mathcal P$. Its degree
in $s$ is less than $d$, so it is $0$. Hence $F^\mu_Q(s\ones+w)=F^{\mu'}_{Q'}(s\ones+w)$, and
Lemma~\ref{lem:reductions}(a) gives the laws at every $T$.
\end{proof}

\subsection{A start without renewal blocks}\label{sec:law-generic}

\begin{definition}[renewal block]\label{def:renewal}
Let $Q$ be irreducible and $\mu$ a start. A set $B$ with $2\le|B|\le d-1$ is a \emph{renewal block}
of $(Q,\mu)$ if $Q[B,B^c]=\ones r^\top$, $Q[B^c,B]=c\rho^\top$ and $\mu_B=m\rho$ for some
$r,c\in[0,\infty)^{B^c}$, a law $\rho$ on $B$ and $m\ge0$. So every state of $B$ leaves $B$ at the
same rates, every jump into $B$ lands with the law $\rho$, and $\mu$ is a multiple of $\rho$ on $B$,
possibly $0$.
\end{definition}

Since $Q$ is irreducible, $c\ne0$, so $\rho$ is determined by $Q$. Put
\[
 M(Q,\mu)=\begin{pmatrix}Q&\ones\\\mu&0\end{pmatrix},
\]
indexed by $[d]\cup\{\partial\}$. A cut of a matrix of order at least $4$ is a partition of its
indices into 2 sets of at least 2 elements such that both off-diagonal blocks have rank at most 1,
as in Loewy~\cite{loewy1986}.

\begin{lemma}[the bordered matrix with the start]\label{lem:bordered-start}
Let $Q$ be irreducible and $\mu$ a start.
\begin{enumerate}[(a)]
\item $\det(\operatorname{diag}(x,t)-M(Q,\mu))=P_Q(x)(t-F^\mu_Q(x))$ with
$P_Q(x)=\det(\operatorname{diag}x-Q)$.
\item For a generator $Q'$ and a start $\mu'$, $O(T)$ has the same law under $(Q',\mu')$ as under
$(Q,\mu)$, for one $T>0$ or equivalently every $T$, if and only if $M(Q',\mu')$ has the principal
minors of $M(Q,\mu)$. Then $\mu'=\mu$.
\item $M(Q,\mu)$ is irreducible, and a set $B\subseteq[d]$ is the side avoiding $\partial$ of a cut of
$M(Q,\mu)$ if and only if $B$ is a renewal block of $(Q,\mu)$.
\end{enumerate}
\end{lemma}

\begin{proof}
(a) For $x$ with large entries the formula of Schur gives $\det(\operatorname{diag}x-Q)\,
(t-\mu(\operatorname{diag}x-Q)^{-1}\ones)$, and both sides are polynomials.

(b) The left side of (a) is multiaffine in $(x,t)$, its coefficients are the principal minors of
$M(Q,\mu)$ up to sign, and its coefficient of $t$ is $P_Q$. So equal principal minors mean
$P_{Q'}=P_Q$ and $F^{\mu'}_{Q'}=F^\mu_Q$. Equal laws give both by Theorem~\ref{thm:one-horizon},
and the converse is Lemma~\ref{lem:reductions}(a). With $m_i=1/x_i$, the expansion
$(\operatorname{diag}x-Q)^{-1}=\sum_{k\ge0}(\operatorname{diag}m\,Q)^k\operatorname{diag}m$ for large
$x$ gives $F^\mu_Q=\mu\cdot m+O(|m|^2)$, so $\mu'=\mu$.

(c) Every $i$ has the edge $i\to\partial$, $\partial$ has edges into the support of $\mu$, and $Q$
is irreducible. Let $J=B^c\cup\{\partial\}$. Both sides have at least 2
elements exactly when $2\le|B|\le d-1$. The column $\partial$ of $M[B,J]$ is $\ones_B$, so $M[B,J]$
has rank at most 1 exactly when $Q[B,B^c]=\ones r^\top$. The rows of $M[J,B]$ are the rows
$Q[y,B]\ge0$ for $y\notin B$, not all $0$, and $\mu_B\ge0$. They span a space of dimension at most
1 exactly when they are multiples $c_y\rho$ and $m\rho$ of one law $\rho$, with $c_y,m\ge0$.
\end{proof}

\begin{theorem}[a start without renewal blocks]\label{thm:generic-start}
Let $Q$ be irreducible, let $\mu$ be a start, and let $Q'$ be a generator and $\mu'$ a start such
that $O(T_0)$ has the same law under $(Q',\mu')$ as under $(Q,\mu)$ for one $T_0>0$. Then
$\mu'=\mu$. If $(Q,\mu)$ has no renewal block, then $Q'=Q$, or $\mu=\pi$ and $Q'=Q^*$. For a fixed
$Q$, every start outside a finite union of proper linear subspaces of $\R^d$ has no renewal block
and differs from $\pi$. So does $\mu=e_y$ when every set $B$ with uniform exit rates and rank-one entries contains
$y$ and has an entry law other than $e_y$.
\end{theorem}

\begin{proof}
Let $M=M(Q,\mu)$ and $M'=M(Q',\mu)$, which have the same principal minors by
Lemma~\ref{lem:bordered-start}(b), and $M$ is irreducible. If $d\ge3$, then $M$ has order at least
$4$ and no cut by Lemma~\ref{lem:bordered-start}(c), and Loewy's
theorem~\cite[Theorem~1]{loewy1986}, with no condition on $M'$, gives an invertible diagonal $E$
with $M'=E^{-1}ME$ or $M'=E^{-1}M^\top E$. If $d\le2$, the theorem of Hartfiel and
Loewy~\cite[Theorem~3]{hartfielloewy1984}, as in the proof of
Corollary~\ref{cor:generator-fiber}(c), gives the same. In the first case the entries $(i,\partial)$
give $1=E_\partial/E_i$, so $E$ is a multiple of $I$ and $Q'=Q$. In the second case they give
$1=\mu_iE_\partial/E_i$, so $\mu_i\ne0$, $E_i=\mu_iE_\partial$ and $q'_{ij}=\mu_jq_{ji}/\mu_i$.
The rows of $Q'$ sum to $0$, so $\mu Q=0$, $\mu=\pi$ and $Q'=Q^*$.

A set $B$ with uniform exit rates and rank-one entries is a renewal block of $(Q,\mu)$ only if
$\mu_B$ lies on the line of $\rho_B$, and this line is a proper subspace of $\R^B$. There are
finitely many such sets, and the starts proportional to $\pi$ form a line as well. For $\mu=e_y$ and such a set $B$,
$\mu_B=e_y$ is not a multiple of $\rho_B$, and $e_y\ne\pi$ since $\pi>0$ and $d\ge2$.
\end{proof}

A renewal block can be lumped to one state. Let $A$ be a renewal block of $(Q,\mu)$ with
$r,c,\rho,m$ as in Definition~\ref{def:renewal} and $\kappa_A=r\cdot\ones>0$, so that the rows $A$
of $Q\ones=0$ give $Q_A\ones=-\kappa_A\ones$. The \emph{lumped pair} $(Q/A,\mu/A)$ lives on
$A^c\cup\{a\}$. The generator $Q/A$ has the off-diagonal entries of $Q$ on $A^c$, the entries
$(Q/A)_{ya}=c_y$ and $(Q/A)_{ay}=r_y$, and zero row sums, and the start is $\mu$ on $A^c$ and $m$
at $a$. Put $g_A(x_A)=\rho(\operatorname{diag}x_A-Q_A)^{-1}\ones$, the transform of one sojourn in
$A$.

\begin{lemma}[lumping]\label{lem:lumping}
For $x\in(0,\infty)^d$ we have $g_A(x_A)>0$, $u_A=1/g_A(x_A)-\kappa_A\ge\min_{i\in A}x_i$ and
$F^\mu_Q(x)=F^{\mu/A}_{Q/A}(x_{A^c},u_A)$.
\end{lemma}

\begin{proof}
The matrix $X_A=\operatorname{diag}x_A-Q_A$ is strictly diagonally dominant with nonpositive
off-diagonal entries, so $X_A^{-1}\ge0$ and $g_A>0$. If $x_A\ge s\ones$, then
$X_A\ones\ge(s+\kappa_A)\ones$, so $g_A\le1/(s+\kappa_A)$. Let $z=(\operatorname{diag}x-Q)^{-1}\ones$
and $z'=z_{A^c}$. The rows $A$ give $X_Az_A=(1+r\cdot z')\ones$, so $\rho\cdot z_A=\zeta$ with
$\zeta=(1+r\cdot z')g_A$. The rows $y\notin A$ involve $z_A$ only through $c_y\,\rho\cdot z_A=c_y\zeta$,
and $(u_A+\kappa_A)\zeta-r\cdot z'=1$. So $(z',\zeta)$ solves
$(\operatorname{diag}(x_{A^c},u_A)-Q/A)(z',\zeta)=\ones$, whose matrix is strictly diagonally
dominant, and $F^\mu_Q=\mu_{A^c}\cdot z'+m\zeta=F^{\mu/A}_{Q/A}(x_{A^c},u_A)$.
\end{proof}

\subsection{The stationary start and memoryless blocks}\label{sec:law-stationary}\label{sec:law-cuts}

\begin{theorem}[the stationary start]\label{thm:stationary-law}
Let $Q$ be irreducible. Then $O(T)$ has the same stationary law under $Q^*$ as under $Q$ for every
$T$.
\end{theorem}

\begin{proof}
By Lemma~\ref{lem:reductions}(a) it suffices that
$F_{Q^*}(x)=\pi D^{-1}(\operatorname{diag}x-Q^\top)^{-1}D\ones=\ones^\top(\operatorname{diag}x-Q^\top)^{-1}\pi^\top=F_Q(x)$.
\end{proof}

\begin{definition}[memoryless block]\label{def:memoryless}
A set $B$ of states with $2\le|B|\le d-1$ is a \emph{memoryless block} of the irreducible $Q$ if
there are $\kappa>0$, a law $\nu$ on $B^c$ and $c\in[0,\infty)^{B^c}$ such that
$q_{iy}=\kappa\nu_y$ and $q_{yi}=c_y\pi_i/\pi(B)$ for all $i\in B$ and $y\notin B$. That is, every
state of $B$ leaves $B$ at the same rate and with the same landing law, and every jump into $B$
lands with the law $\pi_B/\pi(B)$. We write $Q^{[B]}$ for $Q$ with the block $Q_B$ replaced by
$Q^*_B=D_B^{-1}Q_B^\top D_B$, where $D_B$ is the diagonal matrix of $\pi_B$, and the other entries
kept.
\end{definition}

A memoryless block is a strong lumpability. In the terms of Marin and Rossi~\cite{marinrossi2017},
after Buchholz~\cite{buchholz1994}, $B$ is a memoryless block of $Q$ if and only if the partition
into $B$ and the single states outside $B$ is a strong lumpability for $Q$ and for $Q^*$. For $y\notin B$ the first says that
$q_{iy}$, and the second that $\pi_yq_{yi}/\pi_i$, does not depend on $i\in B$.

\begin{theorem}[reversal of a memoryless block]\label{thm:memoryless}
Let $B$ be a memoryless block of the irreducible $Q$. Then $Q^{[B]}$ is an irreducible generator
with stationary law $\pi$, and $F^\mu_{Q^{[B]}}=F^\mu_Q$ for every start $\mu$ with $\mu_B$
proportional to $\pi_B$, possibly $0$. So $O(T)$ has the same law under $(Q,\mu)$ and
$(Q^{[B]},\mu)$ for every $T$, in particular for $\mu=\pi$, and so do the counts of $I+hQ$ and $I+hQ^{[B]}$ for every $N$ and every
admissible $h$. Moreover $Q^{[B]}=Q$ exactly when $Q^*_B=Q_B$, which holds when $|B|=2$, and
$Q^{[B]}=Q^*$ when $|B^c|=1$.
\end{theorem}

\begin{proof}
Put $\rho=\pi_B/\pi(B)$, so $Q[B,B^c]=\kappa\ones\nu^\top$ and $Q[B^c,B]=c\rho$. The flux out of
$B$ equals the flux in, so $\kappa\pi(B)=\sum_y\pi_yc_y$. The columns $B$ of $\pi Q=0$ then give
$\pi_BQ_B=-\kappa\pi_B$, and the rows give $Q_B\ones=-\kappa\ones$. Hence
$Q^*_B\ones=D_B^{-1}Q_B^\top\pi_B^\top=-\kappa\ones$, so the rows of $Q^{[B]}$ sum to $0$, and
$\pi_BQ^*_B=\ones^\top Q_B^\top D_B=-\kappa\pi_B$, which the inflow $\kappa\pi_B$ cancels in the
columns $B$ of $\pi Q^{[B]}$. The other columns do not change. The off-diagonal entries
$\pi_jq_{ji}/\pi_i$ of $Q^*_B$ are nonnegative. The generator $Q^{[B]}$ is irreducible. It has the
same jumps as $Q$ out of $B$, into $B$ and inside $B^c$, and $\kappa>0$ and $c\ne0$ since $Q$ is
irreducible. Take states $s$ and $t$. If $s\in B$, jump out of $B$ first, which is possible from
every state of $B$. A path of $Q$ between two states of $B^c$ stays a path of $Q^{[B]}$ when each of
its sojourns in $B$ is cut to its entry state, since the jumps out of $B$ have the same rates from
every state of $B$. So we reach every state of $B^c$. If $t\in B$, go to a state $y$ with $c_y>0$
and jump to $t$, which is possible since $\rho_t>0$.

Let $\mu_B$ be proportional to $\pi_B$. Then $B$ is a renewal block of $(Q,\mu)$ with entry law
$\rho$, and of $(Q^{[B]},\mu)$ with the same data off $B\times B$, so $Q^{[B]}/B=Q/B$. The sojourn
transform of $Q^{[B]}$ is
$\rho D_B^{-1}(\operatorname{diag}x_B-Q_B^\top)^{-1}D_B\ones=\pi_B(\operatorname{diag}x_B-Q_B)^{-1}\ones/\pi(B)=g_B$,
and Lemma~\ref{lem:lumping} for both gives $F^\mu_{Q^{[B]}}=F^\mu_Q$. The statement on the counts
is Lemma~\ref{lem:reductions}(b).

If $|B|=2$, then $Q_B+\kappa I$ is a generator on 2 states with stationary law $\rho$, hence
reversible. If $B^c=\{y\}$, then $q^*_{iy}=\pi_yc_y/\pi(B)=\kappa=q_{iy}$ and
$q^*_{yi}=\pi_i\kappa/\pi_y=c_y\pi_i/\pi(B)=q_{yi}$, so $Q^{[B]}=Q^*$.
\end{proof}

\begin{remark}[memoryless blocks are sides of cuts]\label{rem:memoryless-cut}
If $|B^c|\ge2$, then $B\cup B^c$ is a cut of $Q$, since $Q[B,B^c]=\kappa\ones\nu^\top$ and
$Q[B^c,B]=c\,\pi_B^\top/\pi(B)$ have rank one.
\end{remark}

\begin{example}[five states]\label{ex:five-state}
Let $d=5$ and $B=\{1,2,3\}$. Inside $B$ put $q_{12}=q_{23}=q_{31}=2$ and $q_{21}=q_{32}=q_{13}=1$.
Each state of $B$ jumps to $4$ at rate $1$ and to $5$ at rate $2$. State $4$ jumps to each state
of $B$ at rate $1$ and to $5$ at rate $1$, and state $5$ jumps to each state of $B$ at rate $2$
and to $4$ at rate $4$. The balance equations give $\pi=(4,4,4,6,3)/21$. So $B$ is a memoryless
block with $\kappa=3$, $\nu=(1/3,2/3)$ and $c=(3,6)$, and since $\pi_B$ is uniform, $Q^{[B]}$ is $Q$
with the block $Q_B$ transposed. It differs from $Q$, since $Q_B$ is not reversible, and from
$Q^*$, since $q^{[B]}_{14}=1$ while $q^*_{14}=\pi_4q_{41}/\pi_1=3/2$. By
Theorem~\ref{thm:memoryless} the stationary laws of $O(T)$ agree for every $T$, and so do the laws
of the counts of $I+hQ$ and $I+hQ^{[B]}$ from $\pi$ for every $N$ and every $h\le1/10$.
\end{example}

The reversals of memoryless blocks are the only obstruction.

\begin{theorem}[the stationary fiber]\label{thm:stationary-fiber}
Let $Q$ and $Q'$ be irreducible generators on $[d]$. Then $O(T)$ has the same law under $Q$ and
under $Q'$, each started from its stationary law, for one $T>0$ or equivalently for every $T$, if
and only if $Q'$ is reached from $Q$ or from $Q^*$ by finitely many reversals of memoryless blocks,
each a memoryless block of the current generator. Theorem~\ref{thm:fiber-count} counts these
generators.
\end{theorem}

The proof puts the stationary transform into a determinant by adding one index. Let $\pi$ be the
stationary law of $Q$, put $e=(\pi_i^{1/2})_{i\in[d]}$ and $K=D_\pi^{1/2}QD_\pi^{-1/2}$, and let
\[
 H(Q)=\begin{pmatrix}K&e\\e^\top&0\end{pmatrix},
\]
a matrix indexed by $[d]\cup\{\partial\}$. Then $Ke=0=e^\top K$, and $H(Q)^\top=H(Q^*)$ since
$K^\top=D_\pi^{1/2}Q^*D_\pi^{-1/2}$. Also $H(Q)=E^{-1}M(Q,\pi)E$ with
$E=\operatorname{diag}(D_\pi^{-1/2},1)$, so by Lemma~\ref{lem:bordered-start} $H(Q)$ is irreducible,
and for 2 generators with the same stationary law, $H(Q)$ and $H(Q')$ have the same principal
minors exactly when the stationary laws of $O(T)$ agree.

\begin{lemma}[the added index selects memoryless blocks]\label{lem:bordered-cuts}
Let $B\subseteq[d]$ with $2\le|B|\le d-1$. Then $B$ is a side of a cut of $H(Q)$ if and only if $B$
is a memoryless block of $Q$. In that case the cut-transpose of $H(Q)$ that transposes the block on
$B$ is $H(Q^{[B]})$, for suitable rank-one factors, and the one that transposes the other side is
$H(Q^{[B]})^\top=H((Q^{[B]})^*)$.
\end{lemma}

\begin{proof}
Put $J=([d]\setminus B)\cup\{\partial\}$. The column $\partial$ of $H(Q)[B,J]$ is $e_B$ and the row
$\partial$ of $H(Q)[J,B]$ is $e_B^\top$, both positive. So the 2 blocks have rank one exactly when
$H(Q)[B,J]=e_Br^\top$ and $H(Q)[J,B]=se_B^\top$ with $r_\partial=s_\partial=1$. For $i\in B$ and
$y\notin B$ this says $q_{iy}=r_ye_y$ and $q_{yi}=s_y\pi_i/e_y$. The first rate does not depend on
$i$, and the second is proportional to $\pi_i$. The total rate $\kappa=\sum_yr_ye_y$ out of $B$ is
positive since $Q$ is irreducible, so $B$ is memoryless with $\nu_y=r_ye_y/\kappa$ and
$c_y=\pi(B)s_y/e_y$. The converse reads the same equations backwards.

We use the cut-transpose $\operatorname{ct}(A,X)$ of~\cite[Definition~2.10]{cggr2024}, defined for
any matrix with a cut $X$. It keeps the block on $X$, transposes the block on the complement, and
replaces the crossing blocks $pq^\top$ on $X\times X^c$ and $uv^\top$ on $X^c\times X$ by $pu^\top$
and $qv^\top$. In the order $(B,J)$ we have
\[
 H(Q)=\begin{pmatrix}K_B&e_Br^\top\\se_B^\top&N\end{pmatrix}.
\]
Take $X=J$, with $p=s$, $q=e_B$, $u=e_B$ and $v=r$. Then $\operatorname{ct}(H(Q),J)$ keeps $N$ and
the crossing blocks $se_B^\top$ and $e_Br^\top$, and replaces $K_B$ by $K_B^\top$. That is
$H(Q^{[B]})$, since $K_B^\top=D_B^{1/2}(D_B^{-1}Q_B^\top D_B)D_B^{-1/2}$. By
Theorem~\ref{thm:memoryless}, $Q^{[B]}$ is an irreducible generator with stationary law $\pi$, so
$H(Q^{[B]})$ has the same border. With $X=B$, $p=e_B$, $q=r$, $u=s$ and $v=e_B$, the same computation
gives $\operatorname{ct}(H(Q),B)=H(Q^{[B]})^\top$. Other factorizations change the result by a
diagonal similarity~\cite[Remark~2.11]{cggr2024}.
\end{proof}

\begin{remark}[without memoryless blocks]\label{rem:no-cggr}
The renewal blocks of $(Q,\pi)$ are the memoryless blocks of $Q$, since $\pi_B=m\rho$ forces
$\rho=\pi_B/\pi(B)$. So if $Q$ has none,
Theorem~\ref{thm:generic-start} with $\mu=\pi$ gives the stationary fiber $\{Q,Q^*\}$, and if
$d\le3$, Corollary~\ref{cor:generator-fiber}(c) applies to $P_{Q'}=P_Q$. These cases need only Loewy,
and Hartfiel and Loewy, and not the theorem of Chatterjee, Ghosh, Gurjar and Raj.
\end{remark}

\begin{proof}[Proof of Theorem~\ref{thm:stationary-fiber}]
The \emph{if} part is Theorems~\ref{thm:memoryless} and~\ref{thm:stationary-law}, applied along the
sequence. Conversely, suppose the stationary laws agree. The first moments give $\pi'=\pi$, and
Theorem~\ref{thm:one-horizon} gives $P_{Q'}=P_Q$ and $F_{Q'}=F_Q$. So $H(Q)$ and $H(Q')$ are
irreducible matrices of order $d+1$ with the same principal minors. By the theorem of Chatterjee,
Ghosh, Gurjar and Raj~\cite[Theorem~1.1]{cggr2024}, there are fewer than $2(d+1)$ cut-transposes
$H(Q)=C_0,C_1,\dots,C_k$, with $C_{s+1}=\operatorname{ct}(C_s,X_s)$ for a cut $X_s$ of $C_s$, such
that $H(Q')$ is diagonally similar to $C_k$ or to $C_k^\top$. A cut there has between 2 and $d-1$
elements on each side~\cite[Definition~2.5]{cggr2024}.

We show by induction that each $C_s$ is diagonally similar to $H(R_s)$ for an irreducible
generator $R_s$ with stationary law $\pi$, obtained from $Q$ by $s$ reversals of memoryless blocks,
each a memoryless block of the current generator, and some time reversals in between. This holds
for $s=0$. A cut of
$C_s$ is a cut of $H(R_s)$, since cuts do not change under diagonal similarity, and a diagonal
similarity of a matrix changes its cut-transpose by a diagonal similarity. Indeed, for
$E=\operatorname{diag}(E_X,E_{X^c})$ a short computation with the blocks gives
$\operatorname{ct}(ECE^{-1},X)=E'\operatorname{ct}(C,X)E'^{-1}$ with
$E'=\operatorname{diag}(E_X,E_{X^c}^{-1})$. Exactly one side of the cut avoids $\partial$, and it has
between 2 and $d-1$ states. By Lemma~\ref{lem:bordered-cuts} it is a memoryless block $B$ of $R_s$,
and $C_{s+1}$ is diagonally similar to $H(R_s^{[B]})$ or to $H((R_s^{[B]})^*)$.

So $H(Q')$ is diagonally similar to $H(R)$ for some $R$ of this form, or to $H(R)^\top=H(R^*)$. If
$H(R_2)=E^{-1}H(R_1)E$ with $E$ diagonal, the entries $(i,\partial)$ give $e_i=E_i^{-1}e_iE_\partial$,
so $E$ is a multiple of the identity and $R_2=R_1$. Hence $Q'$ is $R$ or $R^*$. Finally, a memoryless
block of a generator is a memoryless block of its time reversal, with the rates
$\pi_yc_y/\pi(B)$ out of $B$ and $\pi_i\kappa\nu_y/\pi_y$ into $B$, and $(R^{[B]})^*=(R^*)^{[B]}$,
since both have the entries of $R^*$ off $B\times B$ and those of $R$ on it. So every time reversal
in the sequence can be moved to the front, where 2 of them cancel. This writes $Q'$ as a sequence
of memoryless-block reversals of $Q$ or of $Q^*$, one for each cut-transpose.
\end{proof}

We now count the fiber. Write $\Sigma_{ij}=\pi_iq_{ij}+\pi_jq_{ji}$ and
$\Phi_{ij}=\pi_iq_{ij}-\pi_jq_{ji}$ for $i\ne j$. For a fixed $\pi$, a matrix $R$ with zero row sums
and $\pi R=0$ is determined by its diagonal, $\Sigma^R$ and $\Phi^R$. Write $\rho_B(R)=R^{[B]}$ and
$\rho_{[d]}(R)=R^*$. Then $\rho_B$ keeps the diagonal and $\Sigma$ and changes the sign of
$\Phi_{ij}$ for $i,j\in B$. So the maps $\rho_B$ are involutions and commute with each other, also
for 2 sets that overlap, which means that they meet and neither contains the other. Let
$\mathcal M(Q)$ be the set of memoryless blocks of $Q$.

\begin{lemma}[overlapping blocks]\label{lem:overlap}
Let $B,B'\in\mathcal M(Q)$ overlap, and put $A=B\cap B'$, $X=B\setminus B'$ and $Y=B'\setminus B$.
Then there is $\sigma\ge0$ with $q_{ij}=\sigma\pi_j$ whenever $i$ and $j$ lie in different sets
among $A$, $X$ and $Y$. So $\Phi_{ij}=0$ for such $i,j$, and every union of some of these 3 sets
with at least 2 and at most $d-1$ states is a memoryless block. This includes $B\cup B'$, $B\cap
B'$, $B\setminus B'$ and $B\triangle B'$.
\end{lemma}

\begin{proof}
Put $Z=[d]\setminus(B\cup B')$. There are $r_y,s_y\ge0$ with $q_{iy}=r_y$ and $q_{yi}=s_y\pi_i$ for
$i\in B$ and $y\notin B$, and likewise $r'_y,s'_y$ for $B'$. The entries $q_{az}$ and $q_{za}$ with
$a\in A$ and $z\in Z$ give $r_z=r'_z$ and $s_z=s'_z$. For $x\in X$ and $y\in Y$ we have
$q_{xy}=r_y=s'_x\pi_y$ and $q_{yx}=r'_x=s_y\pi_x$. So $s'_x=\sigma'$ is constant on $X$ and
$s_y=\sigma$ is constant on $Y$. The rates from $X$ to $A\cup Y$ and from $A$ to $Y$ are then
$\sigma'\pi_j$, and those from $Y$ to $A\cup X$ and from $A$ to $X$ are $\sigma\pi_j$. The flux
out of $X$ equals the flux into $X$, which gives $(\sigma'-\sigma)(\pi(A)+\pi(Y))+\Delta=0$ with
$\Delta=\sum_{z\in Z}(r_z-\pi_zs_z)$. The set $Y$ gives $(\sigma-\sigma')(\pi(A)+\pi(X))+\Delta=0$.
Subtracting, $(\sigma'-\sigma)(2\pi(A)+\pi(X)+\pi(Y))=0$, so $\sigma=\sigma'$. Now let $P$ be a union
of some of the sets $A$, $X$ and $Y$. A state of $P$ jumps to a state $j\in(B\cup B')\setminus P$
at rate $\sigma\pi_j$ and to $z\in Z$ at rate $r_z$, and a state outside $P$ jumps to $i\in P$ at a
rate proportional to $\pi_i$. So $P$ is memoryless when $2\le|P|\le d-1$.
\end{proof}

\begin{lemma}[reversal keeps the blocks]\label{lem:invariance}
If $B\in\mathcal M(Q)$, then $\mathcal M(Q^{[B]})=\mathcal M(Q)$. Also $\mathcal M(Q^*)=\mathcal M(Q)$.
\end{lemma}

\begin{proof}
The second claim is in the proof of Theorem~\ref{thm:stationary-fiber}. Let $R=Q^{[B]}$ and
$B'\in\mathcal M(Q)$. Only the entries of $R$ with both indices in $B$ differ from those of $Q$. If
$B'$ is disjoint from $B$ or contains it, the entries $r_{iy}$ and $r_{yi}$ with $i\in B'$ and
$y\notin B'$ are entries of $Q$. If $B'\subsetneq B$ and $y\in B\setminus B'$, then
$r_{iy}=\pi_yq_{yi}/\pi_i=\pi_ys'_y$ and $r_{yi}=\pi_ir'_y/\pi_y$, in the notation of the last
proof. If $B'$ overlaps $B$, the changed entries are $r_{ax}$ and $r_{xa}$ with $a\in B\cap B'$ and
$x\in B\setminus B'$, and Lemma~\ref{lem:overlap} gives $r_{ax}=\pi_xq_{xa}/\pi_a=\sigma\pi_x=q_{ax}$
and $r_{xa}=q_{xa}$. So $\mathcal M(Q)\subseteq\mathcal M(R)$, and since $B\in\mathcal M(R)$ and
$R^{[B]}=Q$, equality holds.
\end{proof}

So every generator of Theorem~\ref{thm:stationary-fiber} has the memoryless blocks of $Q$. Call a member of
$\mathcal M(Q)\cup\{[d]\}$ \emph{strong} if it overlaps no member of $\mathcal M(Q)$. The strong
members form a laminar family $\mathcal L$ that contains $[d]$, so they form a tree under
inclusion with at most $d-1$ members. For $i\ne j$ let $N(i,j)$ be the smallest member of
$\mathcal L$ that contains $i$ and $j$. Call $N\in\mathcal L$ \emph{active} if $\Phi_{ij}\ne0$ for
some $i,j$ with $N(i,j)=N$, and let $k$ be the number of active members. For
$\varepsilon\in\{\pm1\}^{\mathcal L}$ let $Q^\varepsilon$ have the diagonal of $Q$ and the entries
$\pi_iq^\varepsilon_{ij}=(\Sigma_{ij}+\varepsilon_{N(i,j)}\Phi_{ij})/2$ for $i\ne j$.

\begin{theorem}[the size of the stationary fiber]\label{thm:fiber-count}
Let $Q$ be irreducible. The generators $Q'$ of Theorem~\ref{thm:stationary-fiber} are the
$Q^\varepsilon$, and $Q^\varepsilon=Q^\eta$ exactly when $\varepsilon$ and $\eta$ agree on the
active members. So there are exactly $2^k$ of them, and $(\mathbb Z/2)^k$ acts on them simply
transitively by changing the signs of $\Phi$ on the pairs of the active members. Each of them is
obtained from $Q$ or from $Q^*=Q^{-\ones}$ by at most $d-2$ commuting reversals of strong
memoryless blocks. In particular they are only $Q$ and $Q^*$ exactly when $k\le1$.
\end{theorem}

\begin{proof}
Given $\varepsilon$, apply $\rho_N$ for each $N\in\mathcal L$ at which $\varepsilon$ differs from
its value at the parent of $N$, with the value $1$ above $[d]$. For $N(i,j)=P$ the sign of
$\Phi_{ij}$ changes once for each change of $\varepsilon$ along the chain from $P$ up to $[d]$, so
the result is $Q^\varepsilon$. Each $N\ne[d]$ is a memoryless block of the current generator by
Lemma~\ref{lem:invariance}, so $Q^\varepsilon$ is in the fiber by Theorems~\ref{thm:memoryless}
and~\ref{thm:stationary-law}. Starting from $Q^*$ when $\varepsilon_{[d]}=-1$ leaves at most
$|\mathcal L|-1\le d-2$ block reversals.

Conversely, every $B\in\mathcal M(Q)\cup\{[d]\}$ contains disjoint $N_1,\dots,N_r\in\mathcal L$
such that a pair $\{i,j\}$ with $\Phi_{ij}\ne0$ lies in $B$ exactly when it lies in some $N_t$.
This is clear for a strong $B$. Otherwise $B$ overlaps some $B'\in\mathcal M(Q)$, and by
Lemma~\ref{lem:overlap} $B$ is the disjoint union of $B\cap B'$ and $B\setminus B'$, each a
memoryless block or a single state, with $\Phi=0$ on the pairs between them. So induction on $|B|$
applies. Hence $\rho_B(Q^\varepsilon)=\rho_{N_1}\cdots\rho_{N_r}(Q^\varepsilon)$, since both sides
have the diagonal and $\Sigma$ of $Q$ and the same $\Phi$, and this is some $Q^\eta$. By
Theorem~\ref{thm:stationary-fiber} and Lemma~\ref{lem:invariance}, every $Q'$ in the fiber arises from $Q$ or
$Q^*$ by such steps with $B\in\mathcal M(Q)$, so it is some $Q^\varepsilon$. Finally, a pair with
$N(i,j)=N$ and $\Phi_{ij}\ne0$ recovers $\varepsilon_N$ from $Q^\varepsilon$.
\end{proof}

The blocks that carry circulation need not be laminar. If $B$ carries a 3-cycle and $B\cup\{y\}$
and $B\cup\{z\}$ are memoryless, both carry circulation and overlap, and both act as $\rho_B$.

\begin{remark}[modules]\label{rem:modules}
Label the arc from $i$ to $j$ by $q_{ij}/\pi_j$. The memoryless blocks are exactly the sets $M$ with
$2\le|M|\le d-1$ whose outside states see the same labels to and from every $i\in M$, the modules of
these labels, and by Lemma~\ref{lem:overlap} they form, with $[d]$ and the singletons, a partitive
family~\cite{habibpaul2010} whose tree of strong members is $\mathcal L$. A count $2^k$ over a tree
of non-overlapping blocks has the shape of Gallai's theorem on the transitive orientations of a
comparability graph~\cite{gallai1967} and of the reversal of clans of
2-structures~\cite{ehrenfeuchtrozenberg1990a,ehrenfeuchtrozenberg1990b}, a lineage that the
introduction of~\cite{cggr2024} also states.
\end{remark}

\begin{corollary}[when the stationary law sees only time reversal]\label{cor:stationary-cases}
Let $Q$ and $Q'$ be irreducible and give the same stationary law of $O(T)$. Then $Q'$ is $Q$ or
$Q^*$ when $d\le4$, when $Q$ has no cut, and when $Q$ has no memoryless block $B$ with
$|B|\ge3$ and $|B^c|\ge2$. The last holds for all generators with positive rates outside a finite union of proper
linear subspaces, given by equal rates out of 2 states into a third.
\end{corollary}

\begin{proof}
If $|B^c|\ge2$, then $B$ is a side of a cut by Remark~\ref{rem:memoryless-cut}, and when $d\le4$
every memoryless block has $|B|=2$ or $|B^c|=1$. So in each case every member of $\mathcal M(Q)$,
which by Lemma~\ref{lem:invariance} is a memoryless block of every generator $R$ of the fiber, has
$|B|=2$ or $|B^c|=1$, and Theorem~\ref{thm:memoryless} gives $R^{[B]}=R$ or $R^{[B]}=R^*$. So
the fiber is $\{Q,Q^*\}$ by Theorem~\ref{thm:fiber-count}. With positive rates, a memoryless block
needs $q_{iy}=q_{i'y}$ for distinct $i,i'\in B$ and $y\notin B$.
\end{proof}

\subsection{Which starts hide time reversal}\label{sec:law-hiding}

By Theorem~\ref{thm:stationary-law} the stationary start never separates $Q$ from $Q^*$, and by
Theorem~\ref{thm:generic-start} a generic start does. We now find every start that does not. Write
a start as $\mu=\pi\circ f$ with $f\ge0$ and $\pi\cdot f=1$, and put
\[
 \Phi=D_\pi Q-Q^\top D_\pi ,
\]
the matrix of the net probability fluxes $\Phi_{ij}=\pi_iq_{ij}-\pi_jq_{ji}$. It is
antisymmetric, $\Phi\ones=0$, and $\Phi=0$ exactly when $Q$ is reversible. A \emph{circulation
component} of $Q$ is a connected component, with at least 2 states, of the graph on $[d]$ with an
edge $\{i,j\}$ whenever $\Phi_{ij}\ne0$. For a face $F$ we write $A_F=\operatorname{diag}x_F-Q_F$,
a matrix of affine functions of $x\in\R^F$.

\begin{lemma}[the flux form]\label{lem:flux-form}
Let $Q$ be irreducible, $f\in\R^d$ and $\mu=\pi\circ f$. With $u=(\operatorname{diag}x-Q)^{-1}\ones$
and $v=(\operatorname{diag}x-Q)^{-1}f$,
\[
 F^\mu_{Q^*}(x)-F^\mu_Q(x)=u^\top\Phi\,v .
\]
If $F^\mu_{Q^*}=F^\mu_Q$, then $(A_F^{-1}\ones)^\top\Phi_F(A_F^{-1}f_F)=0$ identically in $x\in\R^F$,
for every face $F$, where $\Phi_F$ is the restriction of $\Phi$ to $F\times F$.
\end{lemma}

\begin{proof}
Since $(\operatorname{diag}x-Q^*)^{-1}=D_\pi^{-1}(\operatorname{diag}x-Q^\top)^{-1}D_\pi$, we have
$F^\mu_{Q^*}=f^\top(\operatorname{diag}x-Q^\top)^{-1}\pi^\top=\pi\cdot v$, while
$F^\mu_Q=(\pi\circ f)\cdot u$. Insert $\ones=(\operatorname{diag}x-Q)u$ in the first and
$f=(\operatorname{diag}x-Q)v$ in the second. The terms with $\operatorname{diag}x$ cancel, and
\[
 \pi\cdot v-(\pi\circ f)\cdot u=u^\top D_\pi Qv-v^\top D_\pi Qu=u^\top\Phi v .
\]
For the second claim let $x_j\to\infty$ for $j\notin F$. By the formula of Schur, the block of
$(\operatorname{diag}x-Q)^{-1}$ on $F$ tends to $A_F^{-1}$ and the other blocks tend to $0$. So $u$
and $v$ tend to $A_F^{-1}\ones$ and $A_F^{-1}f_F$ on $F$ and to $0$ off $F$.
\end{proof}

The blocks that matter are memoryless blocks that no net flux crosses. For $X\subseteq[d]$, say that
a state $z\notin X$ is \emph{attached} to $X$ if the rates $q_{iz}$, $i\in X$, are not all equal, or
if $\Phi_{iz}\ne0$ for some $i\in X$.

\begin{lemma}[closed blocks]\label{lem:closure}
\begin{enumerate}[(a)]
\item A set $B$ with $2\le|B|\le d-1$ has no attached state if and only if it is a memoryless block
with $\Phi_{iz}=0$ for all $i\in B$ and $z\notin B$. We call such a block \emph{closed}.
\item Let $|X|\ge2$, put $X_0=X$, and let $X_{k+1}$ be $X_k$ together with the states attached to
$X_k$. The sequence becomes constant at a set $N(X)$. If $X$ lies in a closed block, then $N(X)$ is
the smallest closed block containing $X$, and it lies in every closed block containing $X$.
Otherwise $N(X)=[d]$.
\end{enumerate}
\end{lemma}

\begin{proof}
(a) Let $B$ have no attached state, and write $w_z$ for the common value of $q_{iz}$, $i\in B$, for
$z\notin B$. Since $\Phi_{iz}=0$, $q_{zi}=\pi_iw_z/\pi_z$, which is proportional to $\pi_i$ on $B$.
The total rate $\kappa=\sum_{z\notin B}w_z$ is positive since $Q$ is irreducible, so $B$ is
memoryless with $\nu_z=w_z/\kappa$ and $c_z=\pi(B)w_z/\pi_z$. The converse is the definition.

(b) Let $B$ be a closed block that contains $X_k$, and let $z\notin B$. The rates $q_{iz}$ are equal
for $i\in B$, hence for $i\in X_k$, and $\Phi_{iz}=0$ for $i\in X_k$. So $z$ is not attached to
$X_k$, and every state attached to $X_k$ lies in $B$. By induction every $X_k$ lies in $B$. The
final set $N(X)$ has no attached state, so by (a) it is a closed block or $[d]$.
\end{proof}

A closed block $B$ satisfies $\Phi_{iz}=0$ across its boundary. So if $\Phi$ is supported on
$B\times B$, then $Q^*$ agrees with $Q$ off $B\times B$, and $Q^*=Q^{[B]}$ in the notation of
Definition~\ref{def:memoryless}. The 5 states of Example~\ref{ex:five-state} form a memoryless block
$B=\{1,2,3\}$ that is not closed, since $\Phi_{14}=\pi_1q_{14}-\pi_4q_{41}=-2/21$.

\begin{theorem}[which starts hide time reversal]\label{thm:hiding}
Let $Q$ be irreducible and let $\mu=\pi\circ f$ be a start. Then $O(T)$ has the same law under
$(Q,\mu)$ and under $(Q^*,\mu)$, for one $T>0$ or equivalently for every $T$, if and only if $f$ is
constant on $N(C)$ for every circulation component $C$ of $Q$. Equivalently,
$\E_\pi[f(\xi_T)-f(\xi_0)\mid O(T)]=0$ almost surely for every $T$. In particular, if some
circulation component lies in no closed block, then only the stationary start hides time
reversal.
\end{theorem}

The \emph{if} part is the mechanism of Theorem~\ref{thm:memoryless}. The \emph{only if} part rests on
the next lemma, which we prove by induction along the construction of $N(C)$.

\begin{lemma}[attached states are visible]\label{lem:attached}
Let $C$ be a circulation component. Let $P_0=C$ and $P_j=P_{j-1}\cup\{z_j\}$ for $1\le j\le k$, where
each $z_j$ is attached to $P_{j-1}$. Write $A_j=A_{P_j}$ and $\Phi_j=\Phi_{P_j}$.
\begin{enumerate}[(a)]
\item If $b\in\R^{P_k}$ and $(A_k^{-1}\ones)^\top\Phi_k(A_k^{-1}b)=0$ for all $x$, then $b$ is
constant.
\item For $1\le j\le k$, $(A_j^{-1}\ones)^\top\Phi_j(A_j^{-1}e_{z_j})$ is not identically zero.
\end{enumerate}
\end{lemma}

\begin{proof}
We prove (a) for $k=0$, and then (b) for $j$ and (a) for $k=j$ from (a) for $k=j-1$.

\emph{(a) for $k=0$.} Put $x=t\xi$ with $\xi\in(0,\infty)^C$ fixed and $t\to\infty$. Then
$A_0^{-1}=t^{-1}\operatorname{diag}(\xi)^{-1}+O(t^{-2})$, and the form equals
\[
 t^{-2}\sum_{i,j\in C}\frac{\Phi_{ij}b_j}{\xi_i\xi_j}+O(t^{-3})
 =\frac{t^{-2}}2\sum_{i,j\in C}\frac{\Phi_{ij}(b_j-b_i)}{\xi_i\xi_j}+O(t^{-3}),
\]
by the antisymmetry of $\Phi$. The leading coefficient is a quadratic form in $1/\xi$ that vanishes
on an open set, so $\Phi_{ij}(b_j-b_i)=0$ for all $i,j$. The graph of the edges with $\Phi_{ij}\ne0$
is connected on $C$, so $b$ is constant.

\emph{A Schur complement.} Let $P=P_{j-1}$, $z=z_j$ and $X=P_j$, and write $a=(q_{iz})_{i\in P}$,
$e=(q_{zi})_{i\in P}$, $\phi=(\Phi_{iz})_{i\in P}$ and $s=x_z+q_z$. For $r\in\R^X$ the solution of
$A_Xy'=r$ is $y'=(y,t)$ with
\[
 t=\frac{r_z+e\cdot y}s,\qquad\Bigl(A_P-\frac{ae^\top}s\Bigr)y=r_P+\frac{r_z}s\,a .
\]
For 2 vectors $r,r'$ with solutions $(y,t)$ and $(y',t')$,
\begin{equation}\label{eq:schur-form}
 (A_X^{-1}r)^\top\Phi_X(A_X^{-1}r')=y^\top\Phi_Py'+t'\,y\cdot\phi-t\,y'\cdot\phi .
\end{equation}
As $x_z\to\infty$ with $x_P$ fixed, $s\to\infty$ and $A_P-ae^\top/s\to A_P$.

\emph{(b) for $j$.} Take $r=\ones$ and $r'=e_z$. Then $y=A_P^{-1}\ones+O(1/s)$, $y'=O(1/s)$,
$s\,y'=A_P^{-1}a+O(1/s)$, $t=O(1/s)$ and $s\,t'=1+O(1/s)$. Multiplying~\eqref{eq:schur-form} by $s$
and letting $s\to\infty$, the identity $(A_X^{-1}\ones)^\top\Phi_X(A_X^{-1}e_z)=0$ would give
\begin{equation}\label{eq:attached-limit}
 (A_P^{-1}\ones)^\top\Phi_P(A_P^{-1}a)+(A_P^{-1}\ones)\cdot\phi=0
\end{equation}
for all $x_P$. If $\phi\ne0$, put $x_P=t\xi$ as above. The first term is $O(t^{-2})$ and the second
is $t^{-1}\sum_i\phi_i/\xi_i+O(t^{-2})$, and $\sum_i\phi_i/\xi_i\ne0$ for some $\xi$, a
contradiction. If $\phi=0$, then $z$ is attached through unequal rates, so $a$ is not constant, and
\eqref{eq:attached-limit} contradicts (a) for $k=j-1$.

\emph{(a) for $k=j$.} The form is unchanged when a constant vector is added to $b$, because
$(A_X^{-1}\ones)^\top\Phi_X(A_X^{-1}\ones)=0$. So we may assume $b_z=0$. Take $r=\ones$ and $r'=b$.
Then $y'=A_P^{-1}b_P+O(1/s)$ and $t,t'=O(1/s)$, so~\eqref{eq:schur-form} tends to
$(A_P^{-1}\ones)^\top\Phi_P(A_P^{-1}b_P)$ as $s\to\infty$. By (a) for $k=j-1$, $b_P=\beta\ones$ for
some $\beta$, so $b=\beta\ones-\beta e_z$, and the form of $b$ is $-\beta$ times the form of $e_z$. By
(b) for $j$, $\beta=0$.
\end{proof}

\begin{proof}[Proof of Theorem~\ref{thm:hiding}]
By Lemma~\ref{lem:reductions}(a) and Theorem~\ref{thm:one-horizon}, the laws agree for one $T$ if and
only if they agree for all $T$, and if and only if $F^\mu_{Q^*}=F^\mu_Q$. By
Lemma~\ref{lem:flux-form} this is $u^\top\Phi v=0$ for all $x$. Every edge with $\Phi_{ij}\ne0$ lies
inside one circulation component, so $\Phi=\sum_C\Phi^{(C)}$, where $\Phi^{(C)}$ agrees with $\Phi$ on
$C\times C$ and vanishes elsewhere, and $u^\top\Phi v=\sum_Cu_C^\top\Phi_Cv_C$.

\emph{If.} Fix $C$ and let $B=N(C)$, with $f=c$ on $B$. If $B=[d]$, then $v=cu$ and
$u_C^\top\Phi_Cv_C=c\,u_C^\top\Phi_Cu_C=0$. Otherwise $B$ is a closed block, so
$Q[B,B^c]=\kappa\ones\nu^\top$. The rows $B$ of $(\operatorname{diag}x-Q)u=\ones$ and
$(\operatorname{diag}x-Q)v=f$ give
\[
 u_B=(1+\kappa\,\nu\cdot u_{B^c})A_B^{-1}\ones,\qquad v_B=(c+\kappa\,\nu\cdot v_{B^c})A_B^{-1}\ones .
\]
So $u_C$ and $v_C$ are multiples of the same vector, and $u_C^\top\Phi_Cv_C=0$.

\emph{Only if.} Let $u^\top\Phi v=0$ for all $x$, and fix $C$. By Lemma~\ref{lem:flux-form} with $F=C$
and Lemma~\ref{lem:attached}(a) with $k=0$, $f$ is constant on $C$, say $f=c$ there. List the states
of $N(C)\setminus C$ as $z_1,\dots,z_k$, in the order in which Lemma~\ref{lem:closure}(b) adds them,
in any order within one step. A state added at step $m+1$ is attached to $X_m$, hence to every larger
set, so $z_j$ is attached to $P_{j-1}=C\cup\{z_1,\dots,z_{j-1}\}$. Suppose that $f$ is not constant
on $N(C)$, and let $j$ be the least index with $f_{z_j}\ne c$. Then $f=c$ on $P_{j-1}$, and
$f_{P_j}=c\ones+(f_{z_j}-c)e_{z_j}$. Lemma~\ref{lem:flux-form} with $F=P_j$ gives
$(f_{z_j}-c)(A_j^{-1}\ones)^\top\Phi_j(A_j^{-1}e_{z_j})=0$, which contradicts
Lemma~\ref{lem:attached}(b).

\emph{Conditional form.} By path reversal, $\E^{Q^*}_\mu[\varphi(O(T))]=\E^Q_\pi[f(\xi_T)\varphi(O(T))]$
for bounded $\varphi$, while $\E^Q_\mu[\varphi(O(T))]=\E^Q_\pi[f(\xi_0)\varphi(O(T))]$. The 2 laws agree
for every $\varphi$ exactly when $f(\xi_T)-f(\xi_0)$ has conditional mean $0$ given $O(T)$. Finally,
if $N(C)=[d]$ for some $C$, then $f$ is constant and $\mu=\pi$.
\end{proof}

So the starts that hide time reversal form, after the normalization $\pi\cdot f=1$, the
nonnegative part of the space of functions constant on the blocks $N(C)$. For a generic $Q$ this
space is the constants, and for a reversible $Q$ it is everything.

\subsection{Lumped reversals}\label{sec:law-lumped}

If $A_1$ and $A_2$ are disjoint renewal blocks of $(Q,\mu)$, then $A_2$ is a renewal block of
$(Q/A_1,\mu/A_1)$ with the same $Q_{A_2}$, $\rho$ and $\kappa_{A_2}$. Indeed the rate from $x\in
A_2$ to $a_1$ is $\sum_{i\in A_1}q_{xi}$, which does not depend on $x$, and the rates from $a_1$ into
$A_2$ are the rates $q_{iy}$, $i\in A_1$, $y\in A_2$, which are proportional to $\rho_{A_2}$. So a
family $\mathcal P$ of disjoint renewal blocks can be lumped at once, to $(Q/\mathcal P,\mu/\mathcal P)$
on the atoms, the members of $\mathcal P$ and the other single states, with
$(Q/\mathcal P)_{AB}=\sum_{j\in B}q_{ij}$ for any $i\in A$, and
\begin{equation}\label{eq:lumped-F}
 F^\mu_Q(x)=F^{\mu/\mathcal P}_{Q/\mathcal P}(\tilde x),\qquad \tilde x_{\{j\}}=x_j,\quad
 \tilde x_A=u_A\ (A\in\mathcal P).
\end{equation}
The lumped generator is irreducible, and its stationary law $\pi/\mathcal P$ is $\pi$ summed over
the atoms, since $\sum_{i\in A}\pi_i(Q/\mathcal P)_{AB}=\sum_{i\in A,\,j\in B}\pi_iq_{ij}$.

Conversely, let $\tilde R$ be a generator on the atoms with the diagonal of $Q/\mathcal P$, and let
$\operatorname{un}_{\mathcal P}(\tilde R)$ have the entries $\tilde R_{A(i)A(j)}\rho_{A(j)}(j)$ for $i$
and $j$ in different atoms, with $\rho_{\{j\}}(j)=1$, the block $Q_A$ on each $A\in\mathcal P$, and
the diagonal of $Q$. Its rows sum to $0$, since $\tilde R_{AA}=-\kappa_A$. Each $A\in\mathcal P$ is
a renewal block of $(\operatorname{un}_{\mathcal P}(\tilde R),\mu)$ with the data $Q_A$, $\rho$ and
$\kappa_A$ it has for $Q$, and $\operatorname{un}_{\mathcal P}(\tilde R)/\mathcal P=\tilde R$.

\begin{definition}[lumped reversal]\label{def:lumped-reversal}
Let $\mathcal P$ be a family, possibly empty, of disjoint renewal blocks of $(Q,\mu)$, and put
$\tilde Q=Q/\mathcal P$, $\tilde\mu=\mu/\mathcal P$ and $\tilde\pi=\pi/\mathcal P$. Let $\tilde R$ be
$\tilde Q^{[\tilde B]}$ for a memoryless block $\tilde B$ of $\tilde Q$ with $\tilde\mu_{\tilde B}$
proportional to $\tilde\pi_{\tilde B}$, possibly $0$, or $\tilde Q^*$ when $\tilde\mu=\tilde\pi$.
Then $\operatorname{un}_{\mathcal P}(\tilde R)$ is a \emph{lumped reversal} of $(Q,\mu)$.
\end{definition}

\begin{theorem}[lumped reversals hide]\label{thm:lumped-hides}
Every lumped reversal $R$ of $(Q,\mu)$ is an irreducible generator, and $O(T)$ has the same law
under $(R,\mu)$ as under $(Q,\mu)$ for every $T$. The same holds after finitely many lumped
reversals, each of the current generator.
\end{theorem}

\begin{proof}
By Theorem~\ref{thm:memoryless}, or Theorem~\ref{thm:stationary-law} when $\tilde R=\tilde Q^*$,
$F^{\tilde\mu}_{\tilde R}=F^{\tilde\mu}_{\tilde Q}$. Both choices of $\tilde R$ keep the diagonal of
$\tilde Q$, so $R=\operatorname{un}_{\mathcal P}(\tilde R)$ is a generator with $R/\mathcal P=\tilde R$.
It is irreducible, since $\tilde R$ is (Theorem~\ref{thm:memoryless}), every state of
$A\in\mathcal P$ leaves $A$, and every state of $A$ is reached inside $A$ from the support of
$\rho_A$, as for $Q$. Its blocks $A\in\mathcal P$ are renewal blocks of $(R,\mu)$ with the $g_A$ and
$\kappa_A$ of $Q$. So~\eqref{eq:lumped-F} for $R$ and for $Q$ gives $F^\mu_R=F^\mu_Q$ on
$(0,\infty)^d$, hence as rational functions, and Lemma~\ref{lem:reductions}(a) gives the laws.
\end{proof}

The converse holds, so lumped reversals are the whole fiber from every start. We prove it in
Section~\ref{sec:every-start}.

\begin{theorem}[the fiber from every start]\label{thm:every-start}
Let $Q$ be an irreducible generator and $\mu$ a start. A generator $Q'$ and a start $\mu'$ give
$O(T_0)$ the law it has under $(Q,\mu)$, for one $T_0>0$, if and only if $\mu'=\mu$ and $Q'$ is
reached from $Q$ by finitely many lumped reversals, each a lumped reversal of the current pair.
\end{theorem}

\begin{example}[a lumped reversal that is not a block reversal]\label{ex:lumped}
Let $d=6$, $B_1=\{0,1\}$, $Z=\{2,3\}$, $B_2=B_1\cup Z$, $\rho=(3/4,1/4)$ and $\mu=e_4$. Put
$q_{01}=q_{10}=1$, and for $i\in B_1$ put $q_{i2}=2$, $q_{i3}=1$ and $q_{i4}=q_{i5}=1$. For $z\in Z$
put $q_{zi}=c_z\rho_i$ ($i\in B_1$) with $(c_2,c_3)=(1,2)$, $q_{23}=2$, $q_{32}=1$ and
$q_{z4}=q_{z5}=1$. For $y\in\{4,5\}$ put $q_{yj}=c'_y\,(3/4,1/4,1,1)_j$ ($j\in B_2$) with
$(c'_4,c'_5)=(1,2)$, and $q_{45}=q_{54}=1$. Then $\pi=(171,81,252,252,224,140)/1120$. The renewal
blocks of $(Q,e_4)$ are $B_1$, $B_2$ and $B_2\cup\{5\}$. Their entry laws give the states $0$ and
$1$ the weights $3$ and $1$, while $\pi_0/\pi_1=171/81$, so none is a memoryless block, and
$e_4\ne\pi$. Lump $B_1$ to a state $b$. In $Q/B_1$ the set $\tilde B=\{b,2,3\}$ is left at rate $1$
to each of $4$ and $5$ from each of its states, it is entered from $y$ at the rates
$c'_y(1,1,1)$, and $\tilde\pi_b=\pi_2=\pi_3$. So $\tilde B$ is a memoryless block of $Q/B_1$ with
$\tilde\mu_{\tilde B}=0$, and its rates are $2$ along $b\to2\to3\to b$ and $1$ backwards. The
lumped reversal $Q''$ has $q''_{i2}=1$ and $q''_{i3}=2$ for $i\in B_1$,
$q''_{zi}=(2,1)_z\rho_i$, $q''_{23}=1$, $q''_{32}=2$, and the entries of $Q$ elsewhere. By
Theorem~\ref{thm:lumped-hides} it gives the law of $O(T)$ of $(Q,e_4)$, although $Q''\ne Q$. At
$x=(3,5,7/2,4,9/2,6)$ both transforms are $290948/1250769$.
\end{example}

\subsection{The fiber from every start}\label{sec:every-start}

This subsection proves the \emph{only if} part of Theorem~\ref{thm:every-start}. The \emph{fiber}
of $(Q,\mu)$ is the set of generators $Q'$ such that $O(T)$ has the same law under $(Q',\mu)$ as
under $(Q,\mu)$ for one $T>0$. By Theorem~\ref{thm:one-horizon} and
Lemma~\ref{lem:bordered-start}(b), a pair $(Q',\mu')$ with this property has $\mu'=\mu$, and $Q'$
is in the fiber exactly when $M(Q',\mu)$ has the principal minors of $M(Q,\mu)$. Every $Q'$ in the
fiber is irreducible, $-Q'$ has the principal minors of $-Q$, so $P_{Q'}=P_Q$, and
$F^\mu_{Q'}=F^\mu_Q$. For a set $B$ with \emph{uniform exits and rank-one entries}, that is
$Q[B,B^c]=\ones r^\top$ and $Q[B^c,B]=c\rho^\top$ with $r,c\ge0$ and $\rho$ a law on $B$, we put
$\kappa_B=r\cdot\ones$ and $L_B=Q_B+\kappa_BI$. If $Q$ is irreducible and $2\le|B|\le d-1$, then
$\kappa_B>0$ and $c\ne0$, and the rows $B$ of $Q\ones=0$ give $L_B\ones=0$, so $L_B$ is a generator
on $B$. Every state of $B$ is reached from $\supp\rho$ by $L_B$, since every visit of the chain to
$B$ starts in $\supp\rho$ and moves inside $B$ with the rates of $L_B$. If $A$ is a renewal block of
$(Q,\mu)$, then $Q/A$ is irreducible, since replacing each sojourn of a path of $Q$ in $A$ by the
state $a$ gives a path of $Q/A$.

\begin{lemma}[M-matrices]\label{lem:m-matrix}
Let $R$ be a generator on $[d]$, $S\subseteq[d]$ nonempty and $x\in(0,\infty)^S$. Then
$X=\operatorname{diag}x-R_S$ is invertible, $G=X^{-1}\ge0$ entrywise, $G_{jk}>0$ exactly when
$j=k$ or the graph of $R_S$ has a path from $j$ to $k$, $\det X>0$, and
$G\ones\le(\min_jx_j)^{-1}\ones$.
\end{lemma}

\begin{proof}
Write $X=D-N$ with $D=\operatorname{diag}(x_j-r_{jj})$ and $N\ge0$ the off-diagonal part of $R_S$.
The row $j$ of $N$ sums to at most $-r_{jj}$, so $D^{-1}N$ has row sums at most
$-r_{jj}/(x_j-r_{jj})<1$, and the series $\sum_{m\ge0}(D^{-1}N)^mD^{-1}$ converges to $X^{-1}$. Its
terms are nonnegative, and the entry $(j,k)$ of the term $m$ is positive exactly when there is a
path of length $m$ from $j$ to $k$. Along $x+s\ones$, $s\ge0$, the matrix stays invertible, and
$\det X\to+\infty$ as $s\to\infty$, so $\det X>0$. Finally $X\ones\ge(\min_jx_j)\ones$, since the row
sums of $-R_S$ are nonnegative, and applying $G\ge0$ gives $\ones\ge(\min_jx_j)G\ones$.
\end{proof}

\begin{lemma}[memoryless renewal blocks]\label{lem:memoryless-renewal}
Let $Q$ be irreducible and let $B$, with $2\le|B|\le d-1$, have uniform exits and rank-one entries
with entry law $\rho$. Then $B$ is a memoryless block of $Q$ if and only if $\rho=\pi_B/\pi(B)$, and
if and only if $\rho L_B=0$.
\end{lemma}

\begin{proof}
With $\nu=r/\kappa_B$, Definition~\ref{def:memoryless} asks that $\rho=\pi_B/\pi(B)$. The
columns $B$ of $\pi Q=0$ give $\pi_B(\kappa_BI-L_B)=s\rho$ with $s=\pi_{B^c}\cdot c>0$. The matrix
$\kappa_BI-L_B$ is invertible, since every eigenvalue of the generator $L_B$ has real part at most
$0$. If $\rho L_B=0$, then $\pi_B=(s/\kappa_B)\rho$ solves this equation, so it is the solution. If
$\pi_B=t\rho$, then $\rho L_B=(\kappa_B-s/t)\rho$, and applying both sides to $\ones$ gives
$\kappa_B=s/t$, so $\rho L_B=0$.
\end{proof}

\begin{lemma}[a common renewal block]\label{lem:common-block}
Let $Q'$ be in the fiber of $(Q,\mu)$, and let $A$ be a renewal block of both $(Q,\mu)$ and
$(Q',\mu)$. Then $A$ has the same entry law $\rho$, exit rate $\kappa_A$, transform $g_A$ and
multiple $m$ for $Q$ and $Q'$, and $Q'/A$ is in the fiber of $(Q/A,\mu/A)$. If moreover $A$
contains no smaller renewal block of $(Q,\mu)$, then $Q'_A=Q_A$, or $\rho L_A=0$ and
$Q'_A=D_\rho^{-1}Q_A^\top D_\rho$.
\end{lemma}

\begin{proof}
Write $r,c,\rho,\kappa,m$ for the data of $A$ for $Q$ and $r',c',\rho',\kappa',m'$ for $Q'$. Put
$X_A=\operatorname{diag}x_A-Q_A$ and $X=\operatorname{diag}x_{A^c}-Q_{A^c}$, and $X'_A$, $X'$
likewise for $Q'$. Since $P_{Q'}=P_Q$, the matrices $-Q'$ and $-Q$ have the same principal minors,
so $\det X'_A=\det X_A$ and $\det X'=\det X$. For $x\in(0,\infty)^d$ the formula of Schur on $A$ and
the identity $\det(Y-\gamma cr^\top)=\det Y-\gamma\,r^\top\operatorname{adj}(Y)c$ give
\[
 P_Q(x)=\det X_A\bigl(\det X-g_A\,\mathcal D\bigr),\qquad \mathcal D=r^\top\operatorname{adj}(X)\,c,
\]
since $Q[A^c,A]X_A^{-1}Q[A,A^c]=g_A\,cr^\top$, and likewise $P_{Q'}=\det X_A(\det X-g'_A\mathcal D')$.
Hence $g_A\mathcal D=g'_A\mathcal D'$. Here $\mathcal D>0$ for $x\in(0,\infty)^d$. Indeed some $y$ has
$r_y>0$, and a path of $Q$ from $y$ back to $A$ stays in $A^c$ until it enters $A$ from a state $y'$
with $c_{y'}>0$. So $(X^{-1})_{yy'}>0$ by Lemma~\ref{lem:m-matrix}, $r^\top X^{-1}c>0$ and
$\mathcal D=\det X\,r^\top X^{-1}c>0$. Likewise $\mathcal D'>0$. So $g_A/g'_A=\mathcal D'/\mathcal D$,
where the left side depends only on $x_A$ and the right side only on $x_{A^c}$, and both are one
constant. Since $(tI-Q_A)\ones=(t+\kappa)\ones$, we have $g_A(t\ones)=1/(t+\kappa)$, and likewise
$g'_A(t\ones)=1/(t+\kappa')$, so the constant is $1$ and $\kappa'=\kappa$. Hence $g'_A=g_A$ and
$\mathcal D'=\mathcal D$. With $m_i=1/x_i$, the expansion
$(\operatorname{diag}x_A-Q_A)^{-1}=\sum_{k\ge0}(\operatorname{diag}m\,Q_A)^k\operatorname{diag}m$ gives
$g_A=\rho\cdot m+O(|m|^2)$, so $\rho'=\rho$, and $\mu_A=m\rho=m'\rho$ gives $m'=m$. So $Q/A$ and
$Q'/A$ have the same start $\mu/A$. Expanding along the row $a$,
$\det(\operatorname{diag}(x_{A^c},u)-Q/A)=(u+\kappa)\det X-\mathcal D$, and the same holds for $Q'/A$.
By Lemma~\ref{lem:lumping} for $Q$ and for $Q'$, with $u_A=1/g_A(x_A)-\kappa$, which equals $t$ at
$x_A=t\ones$,
\[
 F^{\mu/A}_{Q'/A}(x_{A^c},t)=F^\mu_{Q'}(x_{A^c},t\ones)=F^\mu_Q(x_{A^c},t\ones)=F^{\mu/A}_{Q/A}(x_{A^c},t)
\]
for all $x_{A^c}\in(0,\infty)^{A^c}$ and $t>0$. So the 2 rational functions agree, and $Q'/A$ is in
the fiber of $(Q/A,\mu/A)$ by Lemma~\ref{lem:reductions}(a).

For the interior, Lemma~\ref{lem:bordered-start}(a) holds with the same proof for every square
matrix in place of $Q$ and every row vector in place of $\mu$. Applied to $Q_A$ and $\rho$, it
gives $\det\bigl(\operatorname{diag}(x_A,t)-N\bigr)=\det X_A\,(t-g_A)$ for
$N=\begin{pmatrix}Q_A&\ones\\\rho&0\end{pmatrix}$, and the same with $Q'_A$ for the analogous $N'$.
So $N$ and $N'$ have the same principal minors, and so have $N+\operatorname{diag}(\kappa\ones_A,0)=
M(L_A,\rho)$ and $M(L'_A,\rho)$, since the principal minors of a matrix determine those of its
diagonal shifts. The matrix $M(L_A,\rho)$ is irreducible, since every state of $A$ is reached from
$\supp\rho$ by $L_A$ and has the edge to $\partial$. Suppose that it has a cut, and let $C\subseteq A$
be the side avoiding $\partial$, so $2\le|C|\le|A|-1$. Then $L_A[C,A\setminus C]=\ones s^\top$, and
the rows of $L_A[A\setminus C,C]$ and $\rho_C$ are multiples of one law $\sigma$ on $C$. So $C$ has
uniform exits for $Q$, since the exits of $A$ to $A^c$ are uniform, its entries from $A^c$ are the
rows $c_y\rho_C$, which are multiples of $\sigma$, and $\mu_C=m\rho_C$ is a multiple of $\sigma$.
So $C$ is a renewal block of $(Q,\mu)$ inside $A$, which is excluded. If $|A|\ge3$, Loewy's
theorem~\cite[Theorem~1]{loewy1986} applies to $M(L_A,\rho)$, and if $|A|=2$ the theorem of
Hartfiel and Loewy~\cite[Theorem~3]{hartfielloewy1984} does, as in the proof of
Theorem~\ref{thm:generic-start}. They give an invertible diagonal $E$ with
$M(L'_A,\rho)=E^{-1}M(L_A,\rho)E$ or $M(L'_A,\rho)=E^{-1}M(L_A,\rho)^\top E$. In the first case the
entries $(i,\partial)$ give $1=E_\partial/E_i$, so $E$ is a multiple of $I$ and $L'_A=L_A$. In the
second case they give $1=\rho_iE_\partial/E_i$, so $\rho_i\ne0$, $E_i=\rho_iE_\partial$ and
$(L'_A)_{ij}=\rho_j(L_A)_{ji}/\rho_i$. The rows of $L'_A$ sum to $0$, so $\rho L_A=0$, and
$L'_A=D_\rho^{-1}L_A^\top D_\rho$. Subtracting $\kappa I$ gives the claim.
\end{proof}

The next step shows that the renewal blocks that contain no smaller one are shared by the whole
fiber. It uses 2 lemmas of Chatterjee, Ghosh, Gurjar and Raj, stated for matrices with nonzero
off-diagonal entries, so we first transfer them. We use their notions: a cut of a matrix of order
$n\ge4$ is a set $X$ of indices with $2\le|X|\le n-2$ such that both off-diagonal blocks between $X$
and its complement have rank at most one~\cite[Definition~2.5]{cggr2024}, so $X$ is a cut exactly
when its complement is, and a cut is minimal if it contains no other cut. Their Lemma~3.5 says the
following. If $A$ and $B$ are $n\times n$ matrices over a field with nonzero off-diagonal entries and
the same principal minors, and $S=\{s_1,s_2\}$ is a cut of $A$, then either $S$ is a cut of $B$, or
for each $i\in\{1,2\}$ the set $X_i=\{t\notin S:A[S\cup\{t\}]$ is diagonally similar to
$B[S\cup\{t\}]\}\cup\{s_i\}$ is a cut of $B$, and $S$ is a cut of $\operatorname{ct}(B,X_i)$. Their
Lemma~3.6 says that under the same hypothesis on $A$ and $B$, a minimal cut of $A$ with more than 2
elements is a cut of $B$.

\begin{lemma}[two lemmas of Chatterjee, Ghosh, Gurjar and Raj, with zero entries]\label{lem:cggr-transfer}
Let $M$ and $M'$ be irreducible real matrices of order $n\ge4$ with the same principal minors.
\begin{enumerate}[(a)]
\item A minimal cut of $M$ with at least 3 elements is a cut of $M'$.
\item If $S=\{s_1,s_2\}$ is a cut of $M$, then $S$ is a cut of $M'$, or $M'$ has the cuts
$P\cup\{s_1\}$ and $P\cup\{s_2\}$ for some set $P$ of indices disjoint from $S$.
\end{enumerate}
\end{lemma}

\begin{proof}
By the theorem of Hartfiel and Loewy~\cite[Theorem~1]{hartfielloewy1984}, every entry of
$\operatorname{adj}(\operatorname{diag}z+M)$ is a nonzero polynomial in $z$, and so is every entry
of $\operatorname{adj}(\operatorname{diag}z+M')$, since $M$ and $M'$ are irreducible. The
determinants $\det(\operatorname{diag}z+M)$ and $\det(\operatorname{diag}z+M')$ are nonzero
polynomials as well, with the coefficient $1$ at $z_1\cdots z_n$. So there is a real diagonal $Z$
such that $M+Z$ and $M'+Z$ are invertible and $\widetilde M=(M+Z)^{\operatorname{adj}}$ and
$\widetilde M'=(M'+Z)^{\operatorname{adj}}$ have nonzero off-diagonal entries. By~\cite[Lemma~2.1]{cggr2024}
they have the same principal minors, and by~\cite[Lemma~2.6]{cggr2024} $\widetilde M$ has exactly
the cuts of $M$ and $\widetilde M'$ exactly those of $M'$. Since the cuts are the same sets, so are
the minimal cuts. Now (a) is~\cite[Lemma~3.6]{cggr2024} for $\widetilde M$ and $\widetilde M'$, and
(b) is the first part of~\cite[Lemma~3.5]{cggr2024} for them, with $X_i=P\cup\{s_i\}$. The second
part of that lemma, on $\operatorname{ct}(\widetilde M',X_i)$, transfers as well
by~\cite[Lemma~2.13]{cggr2024}, but we do not need it.
\end{proof}

The case of 2 states needs more than the cuts: Lemma~\ref{lem:crossing}, proved in
Appendix~\ref{app:two-state}, shows that the second alternative of
Lemma~\ref{lem:cggr-transfer}(b) forces the block anyway.

\begin{proposition}[shared blocks]\label{prop:shared-blocks}
Let $Q'$ be in the fiber of $(Q,\mu)$. Then every renewal block of $(Q,\mu)$ with 2 states, and
every renewal block of $(Q,\mu)$ that contains no smaller one, is a renewal block of $(Q',\mu)$.
\end{proposition}

\begin{proof}
Let $A$ be such a block. Since $|A|\le d-1$, $d\ge3$, and $M=M(Q,\mu)$ and $M'=M(Q',\mu)$ have order
at least $4$. They are irreducible and have the same principal minors
(Lemma~\ref{lem:bordered-start}), and $A$ is a cut of $M$ by Lemma~\ref{lem:bordered-start}(c).
Every cut of $M$ or of $M'$ has exactly one side avoiding $\partial$, and by
Lemma~\ref{lem:bordered-start}(c) the sides avoiding $\partial$ are the renewal blocks of $(Q,\mu)$
and of $(Q',\mu)$.

If $A$ contains no smaller renewal block and $|A|\ge3$, then $A$ is a minimal cut of $M$: a cut
inside $A$ avoids $\partial$, so it is a renewal block inside $A$. By
Lemma~\ref{lem:cggr-transfer}(a) it is a cut of $M'$, hence a renewal block of $(Q',\mu)$.

Let $A=\{1,2\}$. By Lemma~\ref{lem:cggr-transfer}(b), either $A$ is a cut of $M'$, and so a renewal
block of $(Q',\mu)$, or $M'$ has the cuts $X_1=P\cup\{1\}$ and $X_2=P\cup\{2\}$ for a set $P$ of
indices in $([d]\cup\{\partial\})\setminus A$. If $\partial\notin P$, then $X_1$ and $X_2$ avoid
$\partial$, so $P\cup\{1\}$ and $P\cup\{2\}$ are renewal blocks of $(Q',\mu)$, and $P\ne\emptyset$
since $|X_1|\ge2$. If $\partial\in P$, the sides avoiding $\partial$ are the complements $W\cup\{2\}$
and $W\cup\{1\}$ with $W=[d]\setminus(P\cup A)$, and $W\ne\emptyset$ since the complement of $X_1$
has at least 2 elements. In both cases Lemma~\ref{lem:crossing}, with $Y=P$ or $Y=W$, shows that
$A$ is a renewal block of $(Q',\mu)$.
\end{proof}

In particular every renewal block with 2 states is shared by all members of the fiber.

The last tool lifts lumped reversals of a lumped pair. Let $A$ be a renewal block of $(Q,\mu)$ with
the data $\rho$, $\kappa_A$, $m$. For a generator $S$ on $A^c\cup\{a\}$ with $S_{aa}=-\kappa_A$, let
$\operatorname{un}_A(S)$ be the matrix $R$ on $[d]$ with $R_{iy}=S_{ay}$ and $R_{yi}=S_{ya}\rho_i$ for
$i\in A$ and $y\in A^c$, with the off-diagonal entries of $S$ on $A^c\times A^c$, with $R_A=Q_A$, and
with $R_{yy}=S_{yy}$ for $y\in A^c$. Its rows sum to $0$: for $y\in A^c$ since $\rho$ is a law, and
for $i\in A$ since $Q_A\ones=-\kappa_A\ones$ and $\sum_yS_{ay}=\kappa_A$. So $R$ is a generator, $A$ is
a renewal block of $(R,\mu)$ with the data $\rho$, $\kappa_A$, $m$, and $R/A=S$. For $S=Q/A$ it is
$Q$, and for a family $\mathcal P=\{A\}$ it is $\operatorname{un}_{\mathcal P}(S)$ of
Section~\ref{sec:law-lumped}.

\begin{lemma}[lifting]\label{lem:lifting}
Let $S$ be an irreducible generator on $A^c\cup\{a\}$ with $S_{aa}=-\kappa_A$ such that
$R=\operatorname{un}_A(S)$ is irreducible, and let $S'$ be a lumped reversal of $(S,\mu/A)$. Then $\operatorname{un}_A(S')$ is a lumped reversal of
$(R,\mu)$.
\end{lemma}

\begin{proof}
Write $\nu=\mu/A$, let $\mathcal P$ be the family of renewal blocks of $(S,\nu)$ of the lumped
reversal, and $\tilde R$ the generator on the atoms with $S'=\operatorname{un}_{\mathcal P}(\tilde R)$.
For $P\in\mathcal P$ with $a\notin P$ put $\hat P=P$, for $P\ni a$ put $\hat P=(P\setminus\{a\})\cup A$,
and let $\hat{\mathcal P}$ consist of these sets, together with $A$ if no member of $\mathcal P$
contains $a$. The members of $\hat{\mathcal P}$ are pairwise disjoint. We show that they are
renewal blocks of $(R,\mu)$. Recall that $R_{iy}=S_{ay}$ and $R_{yi}=S_{ya}\rho_i$ for $i\in A$ and
$y\in A^c$, that $R$ agrees with $S$ on $A^c\times A^c$ off the diagonal, and that $\mu_A=m\rho$ with
$\nu_a=m$. Let $P\in\mathcal P$ have the entry law $\sigma$ for $S$.

If $a\notin P$, then $P\subseteq A^c$. A state $p\in P$ jumps to $y\in A^c\setminus P$ at the rate
$S_{py}$ and to $i\in A$ at the rate $S_{pa}\rho_i$, and both do not depend on $p$. A state
$y\in A^c\setminus P$ jumps into $P$ with the rates $S_{yP}$, and $i\in A$ with the rates $S_{aP}$,
and both are multiples of $\sigma$, as is $\mu_P=\nu_P$. Also $2\le|P|\le d-1$.

If $a\in P$, then $\hat P$ is left to $y\notin\hat P$ at the rates $S_{py}$ from
$p\in P\setminus\{a\}$ and $S_{ay}$ from $A$, which agree. A state $y\notin\hat P$ jumps to
$j\in P\setminus\{a\}$ at the rate $S_{yj}$ and to $i\in A$ at the rate $S_{ya}\rho_i$, so its rates
into $\hat P$ are $S_{yP}$ spread by $\rho$ over $A$, a multiple of the law $\hat\sigma$ with
$\hat\sigma_j=\sigma_j$ off $A$ and $\hat\sigma_i=\sigma_a\rho_i$ on $A$. Likewise $\mu_{\hat P}$ is
$\nu_P$ spread in the same way, a multiple of $\hat\sigma$. Also $3\le|\hat P|=|P|-1+|A|\le d-1$,
since $|P|\le|A^c|$. Finally $A$ itself is a renewal block of $(R,\mu)$.

The atoms of $\hat{\mathcal P}$ correspond to the atoms of $\mathcal P$, and the rate between 2
atoms is in both cases the sum of the rates into the second atom from any state of the first, so
$R/\hat{\mathcal P}=S/\mathcal P$, with the same start and the same stationary law, $\pi$ summed over
the atoms. So $\tilde R$ is an admissible choice in Definition~\ref{def:lumped-reversal} for
$(R,\mu)$ and $\hat{\mathcal P}$, and $\operatorname{un}_{\hat{\mathcal P}}(\tilde R)$ is a lumped
reversal of $(R,\mu)$. It equals $\operatorname{un}_A(S')$. Both have the blocks of $R$ on the
members of $\hat{\mathcal P}$, since $\operatorname{un}_{\mathcal P}$ keeps the block of $S$ on $P$
and $\operatorname{un}_A$ then restores the block of $R$ on $\hat P$. Both have the diagonal of $R$.
Between different atoms both have the entries $\tilde R_{UV}$ times the entry law of $V$ for $R$,
which is $\hat\sigma$ for $V=\hat P$, $\rho$ for $V=A$ and $1$ for a single state.
\end{proof}

\begin{proof}[Proof of Theorem~\ref{thm:every-start}]
The \emph{if} part and the irreducibility of every generator reached are
Theorem~\ref{thm:lumped-hides}, and $\mu'=\mu$ is Lemma~\ref{lem:bordered-start}(b). For the
\emph{only if} part let $Q'$ be in the fiber of $(Q,\mu)$. We argue by induction on $d$. If
$(Q,\mu)$ has no renewal block, which holds when $d\le2$, then by Theorem~\ref{thm:generic-start}
either $Q'=Q$, or $\mu=\pi$ and $Q'=Q^*$, which is the lumped reversal with $\mathcal P=\emptyset$ and
$\tilde R=Q^*$.

Otherwise let $A$ be a renewal block of $(Q,\mu)$ that contains no smaller one. By
Proposition~\ref{prop:shared-blocks} it is a renewal block of $(Q',\mu)$, so by
Lemma~\ref{lem:common-block} it has the same $\rho$, $\kappa_A$ and $m$ for $Q$ and $Q'$, and $Q'/A$
is in the fiber of $(Q/A,\mu/A)$, a pair on $d-|A|+1<d$ states with an irreducible generator. By
induction there are generators $Q/A=S_0,S_1,\dots,S_k=Q'/A$, each $S_{j+1}$ a lumped reversal of
$(S_j,\mu/A)$. All of them are irreducible and have the diagonal of $Q/A$, since lumped reversals
keep the diagonal, so $S_{j,aa}=-\kappa_A$. By Lemma~\ref{lem:lifting} and induction on $j$, using
Theorem~\ref{thm:lumped-hides} for the irreducibility of $\operatorname{un}_A(S_j)$, each $\operatorname{un}_A(S_{j+1})$ is a lumped
reversal of $(\operatorname{un}_A(S_j),\mu)$, and $\operatorname{un}_A(S_0)=Q$. So
$R=\operatorname{un}_A(Q'/A)$ is reached from $Q$ by lumped reversals. The generators $Q'$ and $R$
have the same lumped generator $Q'/A$ and the same entry law $\rho$ of $A$, so they agree off
$A\times A$, and $R_A=Q_A$. By Lemma~\ref{lem:common-block}, either $Q'_A=Q_A$, and then $Q'=R$, or
$\rho L_A=0$ and $Q'_A=D_\rho^{-1}Q_A^\top D_\rho$. In the second case $L_A=R_A+\kappa_AI$ is also
the interior of $A$ for $R$, so by Lemma~\ref{lem:memoryless-renewal} $A$ is a memoryless block of
$R$ with $\rho=\pi^R_A/\pi^R(A)$, where $\pi^R$ is the stationary law of $R$, and $\mu_A=m\rho$ is
proportional to $\pi^R_A$. So $R^{[A]}$, which replaces $R_A$ by $D_\rho^{-1}R_A^\top D_\rho$, is the
lumped reversal of $(R,\mu)$ with $\mathcal P=\emptyset$ and $\tilde B=A$, and $Q'=R^{[A]}$.
\end{proof}

\section{Computations and open questions}\label{sec:verify}

Every statement is proved in the text. The explicit numbers in Example~\ref{ex:five-state} were
computed in rational arithmetic, and the scripts \texttt{verify\_exactlaw.py} and
\texttt{verify\_reconstruction.py}, which \texttt{verify\_all.py} runs, check them, together with
Corollary~\ref{cor:generator-fiber}(c) on random generators and the
derivative~\eqref{eq:perron-derivative}. The script \texttt{verify\_lumped.py}, also run by
\texttt{verify\_all.py}, checks Example~\ref{ex:lumped} in exact arithmetic: the stationary law,
$Q''\ne Q$, the value of both transforms, and the equality of the 128 principal minors of
$M(Q,e_4)$ and $M(Q'',e_4)$.
Some finite and algebraic steps are also proved in Lean~4 with Mathlib; the code is in the
ancillary files. There the positive null vector of Lemma~\ref{lem:generator-point} is a
hypothesis, since Mathlib has no Perron and Frobenius theorem for irreducible Metzler matrices.
One question remains open: how many visits of
the discrete chain does one stationary count law need to determine the principal minors?

\paragraph{Use of AI in the research.}
We wrote the proofs in this paper. An AI model (Claude Opus~5.5) was used as follows. It suggested
bordering the generator by the start (Lemma~\ref{lem:bordered-start}) and applying Loewy's theorem
to the bordered matrix (Theorem~\ref{thm:generic-start} and Remark~\ref{rem:no-cggr}), and it
proposed a first conjecture for general starts, which was false; we found the counterexample and
proved Theorem~\ref{thm:every-start}. The model was also used to search the literature and check
citations, and to check proofs for errors. The Lean code was generated by Claude Opus~5.5 from our
outlines of the written proofs; it covers the finite and algebraic steps mentioned above and
assumes the positive null vector of Lemma~\ref{lem:generator-point}. We verified all statements
and proofs and take full responsibility for them.

\appendix
\section{Two-state renewal blocks}\label{app:two-state}

This appendix proves Lemma~\ref{lem:crossing}, which Proposition~\ref{prop:shared-blocks} uses for
renewal blocks with 2 states. We first record what $P_{Q'}=P_Q$ says near such a block.

\begin{lemma}[invariants]\label{lem:invariants}
Let $A=\{1,2\}$ have uniform exits and rank-one entries for $Q$, so $q_{1j}=q_{2j}=r_j$ and
$(q_{j1},q_{j2})=c_j(\rho_1,\rho_2)$ for $j\notin A$. Put $\alpha=q_{12}$ and $\beta=q_{21}$, and let
$Q'$ be a generator with $P_{Q'}=P_Q$.
\begin{enumerate}[(a)]
\item $q'_{ii}=q_{ii}$, $q'_{ij}q'_{ji}=q_{ij}q_{ji}$ and
$q'_{ij}q'_{jk}q'_{ki}+q'_{ik}q'_{kj}q'_{ji}=q_{ij}q_{jk}q_{ki}+q_{ik}q_{kj}q_{ji}$ for distinct
$i,j,k$.
\item For $x\in(0,\infty)^{A^c}$ put $G=(\operatorname{diag}x-Q'_{A^c})^{-1}$,
$W=Q'[A,A^c]\,G\,Q'[A^c,A]$, a $2\times2$ matrix, and $w=r^\top(\operatorname{diag}x-Q_{A^c})^{-1}c$.
Then $W_{11}=\rho_1w$, $W_{22}=\rho_2w$ and
$q'_{12}W_{21}+q'_{21}W_{12}-\det W=w(\rho_1\alpha+\rho_2\beta)$.
\end{enumerate}
\end{lemma}

\begin{proof}
(a) The coefficient of $\prod_{k\notin S}x_k$ in $P_Q(x)$ is $\det(-Q_S)$, so $Q$ and $Q'$ have the
same principal minors. Those of order $1$ give the diagonal. With the diagonal equal, those of order
$2$, $q_{ii}q_{jj}-q_{ij}q_{ji}$, give the products. With both equal, the expansion of
$\det(-Q_{\{i,j,k\}})$ leaves the sum of the 2 cyclic products as the only other term.

(b) Fix $x$ and let $x_1,x_2>0$. By Lemma~\ref{lem:m-matrix}, $X=\operatorname{diag}x-Q_{A^c}$ and
$\operatorname{diag}x-Q'_{A^c}$ are invertible, with the same determinant, a principal minor
polynomial of $A^c$. The formula of Schur gives
$P_Q(x_1,x_2,x)=\det X\det\bigl(X_A-Q[A,A^c]X^{-1}Q[A^c,A]\bigr)$, with
$X_A=\operatorname{diag}(x_1,x_2)-Q_A$, and the same for $Q'$, so the 2 determinants of order 2 agree
for all $x_1,x_2>0$. For $Q$ we have $Q[A,A^c]X^{-1}Q[A^c,A]=w\ones\rho^\top$, and
\[
 \det\begin{pmatrix}x_1-q_{11}-w\rho_1&-\alpha-w\rho_2\\-\beta-w\rho_1&x_2-q_{22}-w\rho_2\end{pmatrix}
 =\det X_A-w\rho_2(x_1-q_{11})-w\rho_1(x_2-q_{22})-w(\rho_1\alpha+\rho_2\beta).
\]
For $Q'$ the same expansion gives
$\det X'_A-W_{22}(x_1-q'_{11})-W_{11}(x_2-q'_{22})+\det W-q'_{12}W_{21}-q'_{21}W_{12}$, with
$X'_A=\operatorname{diag}(x_1,x_2)-Q'_A$. Here $\det X'_A=\det X_A$ and $q'_{ii}=q_{ii}$ by (a). After
cancelling $\det X_A$ both sides are affine in $x_1$ and in $x_2$, without a term $x_1x_2$, and they
agree on an open set. The coefficients of $x_1$ and of $x_2$ and the constant terms give the 3
identities.
\end{proof}

\begin{lemma}[the shape of $Q'$]\label{lem:shape}
Let $Q'$ be irreducible, let $Y\subseteq[d]\setminus\{1,2\}$ be nonempty, and let
$B_1=Y\cup\{1\}$ and $B_2=Y\cup\{2\}$ be renewal blocks of $(Q',\mu)$, with entry laws $\sigma_1$,
$\sigma_2$ and multiples $m_1$, $m_2$. Put $U=Y\cup\{1,2\}$ and $K=[d]\setminus U$, and write
$q'_{iY}$ for the row vector $(q'_{iy})_{y\in Y}$ and $|q'_{iY}|$ for its sum.
\begin{enumerate}[(a)]
\item There are $\varepsilon\in[0,\infty)^K$ and $a,b\ge0$ with $q'_{ik}=\varepsilon_k$ for $i\in U$
and $k\in K$, and $q'_{y2}=q'_{12}=a$ and $q'_{y1}=q'_{21}=b$ for $y\in Y$. So
$q'_{11}=-(a+|q'_{1Y}|+|\varepsilon|)$ and $q'_{22}=-(b+|q'_{2Y}|+|\varepsilon|)$.
\item $(b,q'_{2Y})=c^1_2\sigma_1$, $(q'_{k1},q'_{kY})=c^1_k\sigma_1$, $(a,q'_{1Y})=c^2_1\sigma_2$ and
$(q'_{k2},q'_{kY})=c^2_k\sigma_2$ for $k\in K$, with all $c^i_j\ge0$, and $\mu_{B_i}=m_i\sigma_i$.
\item If $b=0$, then $q'_{2Y}=0$. If $a=0$, then $q'_{1Y}=0$.
\item If $b>0$ and $q'_{2Y}=0$, then no state of $K\cup\{2\}$ jumps into $Y$. If $a>0$ and
$q'_{1Y}=0$, then no state of $K\cup\{1\}$ jumps into $Y$.
\item If $q'_{1Y}=q'_{2Y}\ne0$, or if $a=b=0$, then $\{1,2\}$ is a renewal block of $(Q',\mu)$.
\end{enumerate}
\end{lemma}

\begin{proof}
(a) The uniform exits of $B_1$ say that the rows $i\in B_1$ agree on $K\cup\{2\}$, and those of
$B_2$ that the rows $i\in B_2$ agree on $K\cup\{1\}$. Since $Y\ne\emptyset$, the rows $1$, $2$ and
$y\in Y$ agree on $K$, the rows $1$ and $y$ agree at $2$, and the rows $2$ and $y$ agree at $1$. The
diagonal makes the rows sum to $0$.

(b) This is the condition on the entries. The rows of $Q'[B_1^c,B_1]$ are the rows $2$ and $k\in K$
restricted to $B_1$, and likewise for $B_2$.

(c) Let $b=0$ and suppose that $q'_{2Y}\ne0$. By (b), $c^1_2>0$ and $\sigma_1(1)=0$, so no state of
$K\cup\{2\}$ jumps to $1$, and by (a) no $y\in Y$ does. So no state jumps to $1$, which contradicts
irreducibility. The second claim is the same with $1$ and $2$ exchanged.

(d) By (b), $c^1_2\sigma_1=(b,0)$ with $b>0$, so $\sigma_1$ is the point mass at $1$, and by (b)
again $q'_{kY}=0$ for $k\in K$. The second claim is the same with $1$ and $2$ exchanged.

(e) First let $v=q'_{1Y}=q'_{2Y}\ne0$. By (c), $a>0$ and $b>0$. By (b),
$\sigma_1=(b,v)/(b+|v|)$ and $\sigma_2=(a,v)/(a+|v|)$. For $k\in K$, (b) gives
$q'_{kY}=c^1_kv/(b+|v|)=c^2_kv/(a+|v|)$, so $c^1_k/(b+|v|)=c^2_k/(a+|v|)=\lambda_k$ and
$(q'_{k1},q'_{k2})=\lambda_k(b,a)$. For $y\in Y$, $(q'_{y1},q'_{y2})=(b,a)$ by (a). Likewise
$\mu_Y=m_1v/(b+|v|)=m_2v/(a+|v|)$ gives $m_1/(b+|v|)=m_2/(a+|v|)=\lambda$ and $\mu_{\{1,2\}}=\lambda(b,a)$.
The rows $1$ and $2$ agree off $\{1,2\}$, with $v$ on $Y$ and $\varepsilon$ on $K$. So $\{1,2\}$ has
uniform exits, its entries are multiples of the law $(b,a)/(a+b)$, and so is $\mu_{\{1,2\}}$, and
$2\le|\{1,2\}|\le d-1$ since $Y\ne\emptyset$.

Now let $a=b=0$. By (c), $q'_{1Y}=q'_{2Y}=0$, so the rows $1$ and $2$ vanish on $U$ off the
diagonal and are $\varepsilon$ on $K$, which are uniform exits. Every jump into $Y$ from outside $Y$
comes from $K$, and by irreducibility there is one, from some $k_0$. By (b),
$c^1_{k_0}\sigma_1|_Y=c^2_{k_0}\sigma_2|_Y\ne0$, so $\sigma_2|_Y=\lambda\sigma_1|_Y$ with
$\lambda=c^1_{k_0}/c^2_{k_0}>0$. For every $k\in K$, $c^1_k\sigma_1|_Y=\lambda c^2_k\sigma_1|_Y$ gives
$c^1_k=\lambda c^2_k$, so $(q'_{k1},q'_{k2})=c^2_k\theta$ with $\theta=(\lambda\sigma_1(1),\sigma_2(2))$,
and the rows $y\in Y$ vanish on $\{1,2\}$ by (a). From $\mu_Y=m_1\sigma_1|_Y=m_2\lambda\sigma_1|_Y$ we
get $m_1=\lambda m_2$ and $\mu_{\{1,2\}}=m_2\theta$. Finally $\theta\ne0$, since otherwise no state
outside $\{1,2\}$ jumps into it, which contradicts irreducibility. So $\{1,2\}$ is a renewal block
with entry law $\theta/(\theta\cdot\ones)$.
\end{proof}

\begin{lemma}[crossing blocks]\label{lem:crossing}
Let $Q$ and $Q'$ be irreducible with $P_{Q'}=P_Q$, and let $A=\{1,2\}$ be a renewal block of
$(Q,\mu)$. Let $Y\subseteq[d]\setminus A$ be nonempty, and suppose that $Y\cup\{1\}$ and
$Y\cup\{2\}$ are renewal blocks of $(Q',\mu)$. Then $A$ is a renewal block of $(Q',\mu)$.
\end{lemma}

\begin{proof}
Use the notation of Lemmas~\ref{lem:invariants} and~\ref{lem:shape}, and put $R_y=r_yc_y\ge0$ for
$y\in Y$ and $s=|q'_{1Y}|$. For $y\in Y$, Lemma~\ref{lem:invariants}(a) on the pairs $\{1,2\}$,
$\{1,y\}$ and $\{2,y\}$, on the triple $\{1,2,y\}$ and on the diagonal entries $1$ and $2$, with
Lemma~\ref{lem:shape}(a), gives
\begin{align}
 &\alpha\beta=ab, \label{eq:I1}\\
 &R_y\rho_1=b\,q'_{1y},\qquad R_y\rho_2=a\,q'_{2y}, \label{eq:I2}\\
 &R_y(\rho_1\alpha+\rho_2\beta)=ab\,(q'_{1y}+q'_{2y}), \label{eq:I3}\\
 &(a+|q'_{1Y}|)-(b+|q'_{2Y}|)=\alpha-\beta. \label{eq:I4}
\end{align}
For~\eqref{eq:I2}, $q_{1y}q_{y1}=r_yc_y\rho_1$ and $q'_{1y}q'_{y1}=q'_{1y}b$, and likewise for $2$.
For~\eqref{eq:I3}, the cyclic products are $\alpha r_yc_y\rho_1+r_yc_y\rho_2\beta$ for $Q$ and
$aq'_{2y}b+q'_{1y}ab$ for $Q'$. For~\eqref{eq:I4}, $q_{11}-q_{22}=\beta-\alpha$, since the rows $1$
and $2$ of $Q$ have the same exit total $r\cdot\ones$, and $q'_{11}-q'_{22}$ is read off
Lemma~\ref{lem:shape}(a).

\emph{Case $a>0$ and $b>0$.} By~\eqref{eq:I2}, $q'_{1y}=\rho_1R_y/b$ and $q'_{2y}=\rho_2R_y/a$. If
$\Sigma=\sum_{y\in Y}R_y$ were $0$, then $q'_{1Y}=q'_{2Y}=0$, and by both halves of
Lemma~\ref{lem:shape}(d) no state outside $Y$ would jump into $Y$, which contradicts
irreducibility. So $R_{y_0}>0$ for some $y_0$, and~\eqref{eq:I3} at $y_0$, divided by $R_{y_0}$,
reads $\rho_1\alpha+\rho_2\beta=ab(\rho_1/b+\rho_2/a)=\rho_1a+\rho_2b$. With~\eqref{eq:I1}, the
numbers $\rho_1\alpha$ and $\rho_2\beta$ have the sum $\rho_1a+\rho_2b$ and the product
$\rho_1\rho_2ab$, so they are $\rho_1a$ and $\rho_2b$ in some order. Equation~\eqref{eq:I4} reads
\begin{equation}\label{eq:I4b}
 a-b+\Sigma\Bigl(\frac{\rho_1}{b}-\frac{\rho_2}{a}\Bigr)=\alpha-\beta .
\end{equation}
If $\rho_2=0$, then $\rho_1=1$ and $\alpha\in\{a,0\}$. Since $\alpha\beta=ab>0$, $\alpha=a$ and
$\beta=b$, and~\eqref{eq:I4b} gives $\Sigma/b=0$, a contradiction. The case $\rho_1=0$ is the same
with $1$ and $2$ exchanged. If $\rho_1\rho_2>0$ and $(\alpha,\beta)=(a,b)$, then~\eqref{eq:I4b}
gives $\Sigma(\rho_1a-\rho_2b)/(ab)=0$, so $\rho_1a=\rho_2b$. If $\rho_1\rho_2>0$ and
$(\rho_1\alpha,\rho_2\beta)=(\rho_2b,\rho_1a)$, then, since $\rho_1+\rho_2=1$,
$\alpha-\beta-a+b=b/\rho_1-a/\rho_2=(\rho_2b-\rho_1a)/(\rho_1\rho_2)$, and~\eqref{eq:I4b} gives
$(\rho_1a-\rho_2b)\bigl(\Sigma/(ab)+1/(\rho_1\rho_2)\bigr)=0$, so again $\rho_1a=\rho_2b$. Hence
$q'_{1y}=\rho_1R_y/b=\rho_2R_y/a=q'_{2y}$ for all $y\in Y$, and $q'_{1Y}=q'_{2Y}\ne0$ since
$\Sigma>0$. Lemma~\ref{lem:shape}(e) gives the claim.

\emph{Case $b=0<a$.} By Lemma~\ref{lem:shape}(c), $q'_{2Y}=0$, so the row $2$ of $Q'$ vanishes on
$U\setminus\{2\}$, and by Lemma~\ref{lem:shape}(a) no state of $U\setminus\{1\}$ jumps to $1$.
By~\eqref{eq:I1}, $\alpha\beta=0$, and by~\eqref{eq:I4}, $\alpha-\beta=a+s>0$, so $\beta=0$ and
$\alpha=a+s$. Now take $G$ and $W$ as in Lemma~\ref{lem:invariants}(b), and split $A^c=Y\cup K$.
The rows of $Q'[A,A^c]$ are $(q'_{1Y},\varepsilon)$ and $(0,\varepsilon)$, and its columns
$Q'[A^c,A]$ are $(0,q'_{K1})$ and $(a\ones_Y,q'_{K2})$. Hence
\[
 W_{11}=W_{21}+\Delta,\quad W_{12}=W_{22}+\Gamma,\qquad
 \Delta=q'_{1Y}G_{YK}q'_{K1},\quad \Gamma=q'_{1Y}\bigl(aG_{YY}\ones+G_{YK}q'_{K2}\bigr),
\]
and $\det W=\Delta W_{22}-\Gamma W_{21}$. The last identity of Lemma~\ref{lem:invariants}(b), with
$q'_{12}=a$, $q'_{21}=0$, $\beta=0$ and $w\rho_1=W_{11}$, reads $aW_{21}-\det W=(a+s)W_{11}$, that is
\begin{equation}\label{eq:key}
 (\Gamma-s)\,W_{21}=\Delta\,(W_{22}+a+s).
\end{equation}
Suppose that $s>0$, and take $x=t\ones$ with $t>a+\max_{k\in K}q'_{k2}$, the maximum being $0$ if
$K=\emptyset$. By Lemma~\ref{lem:m-matrix}, $G\ge0$ and $G\ones\le t^{-1}\ones$, so
$\Gamma\le s(a+\max_kq'_{k2})/t<s$, while $W_{21},W_{22},\Delta\ge0$. So the left side
of~\eqref{eq:key} is at most $0$ and the right side at least $0$, and both vanish. Hence
$W_{21}=0$ and, since $W_{22}+a+s>0$, $\Delta=0$, so $W_{11}=0$. But $Q'$ is irreducible, so there is
a path $1=i_0\to i_1\to\dots\to i_n=1$ with $n\ge2$. Its last step comes from $K$, since only states
of $K$ jump to $1$. Let $p$ be the last index below $n$ with $i_p\in A$, so $p<n-1$ and
$i_{p+1},\dots,i_{n-1}\in A^c$. By Lemma~\ref{lem:m-matrix},
$q'_{i_pi_{p+1}}G_{i_{p+1}i_{n-1}}q'_{i_{n-1}1}>0$, and this is one of the nonnegative terms of
$W_{i_p1}$. So $W_{11}+W_{21}>0$, a contradiction. Hence $s=0$, that is $q'_{1Y}=0=q'_{2Y}$ with
$a>0$. By Lemma~\ref{lem:shape}(d), no state of $K\cup\{1\}$ jumps into $Y$, and neither does $2$,
which contradicts irreducibility. So this case does not occur.

\emph{Case $a=0<b$.} This is the previous case with $1$ and $2$ exchanged, and it does not occur.

\emph{Case $a=b=0$.} Lemma~\ref{lem:shape}(e) gives the claim.
\end{proof}

\section*{Data and code availability}

The Python scripts and the Lean library of Section~\ref{sec:verify} are ancillary files of the
arXiv version. They are also archived at Zenodo,
\href{https://doi.org/10.5281/zenodo.24099074}{doi:10.5281/zenodo.24099074}.

The scripts need Python~3 with \texttt{numpy}, \texttt{scipy}, \texttt{sympy} and \texttt{mpmath}.
Running \texttt{python3 verify\_all.py} repeats the checks of Section~\ref{sec:verify}. It saves the
output of each script together with the sha256 of that script.

The Lean library is in the folder \texttt{lean}. It uses Lean~4.33.1 and Mathlib~v4.33.1, and the
files \texttt{lean-toolchain} and \texttt{lake-manifest.json} fix these versions. To build it, we
run \texttt{lake exe cache get} and then \texttt{lake build PaperC}. The file \texttt{MAP.md} lists
the statements of this paper that each Lean file proves. In the folder \texttt{lean}, the command
\texttt{python3 verify\_lean.py PaperC --workspace .} builds the library and checks that each
theorem uses only the axioms \texttt{propext}, \texttt{Classical.choice} and \texttt{Quot.sound}.

\section*{Funding}

This research did not receive any specific grant from funding agencies in the
public, commercial, or not-for-profit sectors.

\section*{Declaration of competing interest}

The author declares that he has no known competing financial interests or
personal relationships that could have appeared to influence the work reported
in this paper.

\section*{Declaration of generative AI and AI-assisted technologies in the manuscript preparation process}

During the preparation of this work the author used Claude Opus~5.5 to proofread, to check
notation for consistency, to revise the language of text written by the author, to reorganize
sections and move material between them, and to review drafts. These reviews suggested the present
framing of the main theorems and informed the phrasing of the abstract. The use of the tool in the
research is described in Section~\ref{sec:verify}. After using this tool, the
author reviewed and edited the content and takes full responsibility for the content of the
published article.

\section*{Acknowledgments}

This work grew out of Problem~131 of Peter Kagey's Open Problem
Collection~\cite{kagey131}. We thank him for posing it.


{\small
\begin{thebibliography}{99}
\setlength{\itemsep}{0pt}

\bibitem{alahmadieh2023}
A.~Al Ahmadieh, The fiber of the principal minor map, preprint, arXiv:2309.00806v2, 2023.
\url{https://arxiv.org/abs/2309.00806v2}.

\bibitem{bermanplemmons1994}
A.~Berman, R.~J.~Plemmons, Nonnegative Matrices in the Mathematical Sciences, SIAM,
Philadelphia, 1994. \url{https://doi.org/10.1137/1.9781611971262}.

\bibitem{bfg2014}
L.~Bertini, A.~Faggionato, D.~Gabrielli, From level 2.5 to level 2 large deviations for
continuous time Markov chains, Markov Process. Related Fields 20 (2014) 545--562.
\url{https://arxiv.org/abs/1212.6908}.

\bibitem{bfg2015}
L.~Bertini, A.~Faggionato, D.~Gabrielli, Large deviations of the empirical flow for continuous
time Markov chains, Ann. Inst. Henri Poincar\'e Probab. Stat. 51 (2015) 867--900.
\url{https://doi.org/10.1214/14-AIHP601}.

\bibitem{brunoyangpearson2005}
W.~J.~Bruno, J.~Yang, J.~E.~Pearson, Using independent open-to-closed transitions to simplify
aggregated Markov models of ion channel gating kinetics, Proc. Natl. Acad. Sci. USA 102 (2005)
6326--6331. \url{https://doi.org/10.1073/pnas.0409110102}.

\bibitem{bhk2007}
D.~Brydges, R.~van der Hofstad, W.~K\"onig, Joint density for the local times of continuous-time
Markov chains, Ann. Probab. 35 (2007) 1307--1332. \url{https://doi.org/10.1214/009171906000001024}.

\bibitem{buchholz1994}
P.~Buchholz, Exact and ordinary lumpability in finite Markov chains, J. Appl. Probab. 31 (1994)
59--75. \url{https://doi.org/10.2307/3215235}.

\bibitem{cggr2024}
A.~Chatterjee, S.~Ghosh, R.~Gurjar, R.~Raj, Characterizing and testing principal minor
equivalence of matrices, in: Proc. 57th Annual ACM Symposium on Theory of Computing (STOC 2025),
ACM, New York, 2025, pp.~1067--1078. \url{https://doi.org/10.1145/3717823.3718146}. Full version:
preprint, arXiv:2410.01961v1, 2024, \url{https://arxiv.org/abs/2410.01961v1}. Results, including
Lemmas~3.5 and~3.6, are cited with the numbering of the arXiv version.

\bibitem{cohen1981}
J.~E.~Cohen, Convexity of the dominant eigenvalue of an essentially nonnegative matrix, Proc.
Amer. Math. Soc. 81 (1981) 657--658. \url{https://doi.org/10.1090/S0002-9939-1981-0601750-2}.

\bibitem{ddz2024}
J.-F.~Delmas, D.~Dronnier, P.-A.~Zitt, Transformations preserving the effective spectral radius
of a matrix, Linear Multilinear Algebra 73 (2025) 2354--2360.
\url{https://doi.org/10.1080/03081087.2024.2447523}.

\bibitem{dv1975}
M.~D.~Donsker, S.~R.~S.~Varadhan, On a variational formula for the principal eigenvalue for
operators with maximum principle, Proc. Natl. Acad. Sci. USA 72 (1975) 780--783.
\url{https://doi.org/10.1073/pnas.72.3.780}.

\bibitem{dv1976}
M.~D.~Donsker, S.~R.~S.~Varadhan, On the principal eigenvalue of second-order elliptic
differential operators, Comm. Pure Appl. Math. 29 (1976) 595--621.
\url{https://doi.org/10.1002/cpa.3160290606}.

\bibitem{ehrenfeuchtrozenberg1990a}
A.~Ehrenfeucht, G.~Rozenberg, Theory of 2-structures, part I: clans, basic subclasses, and
morphisms, Theoret. Comput. Sci. 70 (1990) 277--303.
\url{https://doi.org/10.1016/0304-3975(90)90129-6}.

\bibitem{ehrenfeuchtrozenberg1990b}
A.~Ehrenfeucht, G.~Rozenberg, Theory of 2-structures, part II: representation through labeled
tree families, Theoret. Comput. Sci. 70 (1990) 305--342.
\url{https://doi.org/10.1016/0304-3975(90)90130-A}.

\bibitem{fredkinrice1986}
D.~R.~Fredkin, J.~A.~Rice, On aggregated Markov processes, J. Appl. Probab. 23 (1986) 208--214.
\url{https://doi.org/10.2307/3214130}.

\bibitem{friedland1981}
S.~Friedland, Convex spectral functions, Linear Multilinear Algebra 9 (1981) 299--316.
\url{https://doi.org/10.1080/03081088108817381}.

\bibitem{gallai1967}
T.~Gallai, Transitiv orientierbare Graphen, Acta Math. Acad. Sci. Hungar. 18 (1967) 25--66.
\url{https://doi.org/10.1007/BF02020961}.

\bibitem{good1961}
I.~J.~Good, The frequency count of a Markov chain and the transition to continuous time, Ann.
Math. Statist. 32 (1961) 41--48. \url{https://doi.org/10.1214/aoms/1177705138}.

\bibitem{habibpaul2010}
M.~Habib, C.~Paul, A survey of the algorithmic aspects of modular decomposition, Comput. Sci. Rev.
4 (2010) 41--59. \url{https://doi.org/10.1016/j.cosrev.2010.01.001}.

\bibitem{hartfielloewy1984}
D.~J.~Hartfiel, R.~Loewy, On matrices having equal corresponding principal minors, Linear
Algebra Appl. 58 (1984) 147--167. \url{https://doi.org/10.1016/0024-3795(84)90209-X}.

\bibitem{kagey131}
P.~Kagey, Problem 131, Open Problem Collection. \url{https://peterkagey.com/problems/131/},
2026 (accessed October 2026).

\bibitem{kl2024}
S.~Kim, C.~Landim, $\Gamma$-expansion of the measure-current large deviations rate functional of
non-reversible finite-state Markov chains, Stochastic Process. Appl. 189 (2025) 104733.
\url{https://doi.org/10.1016/j.spa.2025.104733}. Section~5 is cited with the numbering of
arXiv:2412.13515v2.

\bibitem{larget1998}
B.~Larget, A canonical representation for aggregated Markov processes, J. Appl. Probab. 35 (1998)
313--324. \url{https://doi.org/10.1239/jap/1032192850}.

\bibitem{loewy1986}
R.~Loewy, Principal minors and diagonal similarity of matrices, Linear Algebra Appl. 78 (1986)
23--64. \url{https://doi.org/10.1016/0024-3795(86)90015-7}.

\bibitem{marinrossi2017}
A.~Marin, S.~Rossi, On the relations between Markov chain lumpability and reversibility, Acta
Inform. 54 (2017) 447--485. \url{https://doi.org/10.1007/s00236-016-0266-1}.

\bibitem{pemmasaniRS}
A.~Pemmasani, Occupation counts of rarely switching Markov chains: stratified local laws, support
selection and the boundary layer, preprint, arXiv:2610.09258v2, 2026.
\url{https://arxiv.org/abs/2610.09258}.

\bibitem{pemmasaniC3}
A.~Pemmasani, What occupation statistics determine about a diffusion, in preparation.

\bibitem{ryden1996}
T.~Ryd\'en, On identifiability and order of continuous-time aggregated Markov chains,
Markov-modulated Poisson processes, and phase-type distributions, J. Appl. Probab. 33 (1996)
640--653. \url{https://doi.org/10.2307/3215346}.

\bibitem{xzdy2024}
X.~Xiang, J.~Zhou, Y.~Deng, X.~Yang, Identifying the generator matrix of a stationary Markov chain
using partially observable data, Chaos 34 (2024) 023132. \url{https://doi.org/10.1063/5.0156458}.

\end{thebibliography}
}
\end{document}